\input amstex
\documentstyle{amsppt}
\magnification=1200
\hcorrection{.25in}
\advance\vsize-.75in
\input epsf
\NoBlackBoxes
\topmatter
\title
3-manifolds as Viewed From The Curve Complex
\endtitle
\author
John Hempel
\endauthor
\affil
Rice University, Houston, TX and MSRI, Berkeley, CA
\endaffil
\address
Rice University, Houston, Texas
\endaddress
\email
hempel\@math.rice.edu
\endemail
\thanks
Research at MSRI is supported in part by NSF grant DMS-9022140
\endthanks
\abstract
A Heegaard diagram for a 3-manifold is regarded as a pair of simplexes
in the complex of curves on a surface and a Heegaard splitting as a
pair of subcomplexes generated by the equivalent diagrams. We relate
 geometric and combinatorial properties of these subcomplexes with
topological properties of the manifold and/or the associated
splitting. For example we show that for any splitting of a 3-manifold
which is Seifert fibered or which contains an essential torus the
subcomplexes are at a distance at most two apart in the simplicial
distance on the curve complex; whereas there are splittings in which
the subcomplexes are arbitrarily far apart. We also give obstructions,
computable from a given diagram, to being Seifert fibered or to
containing an essential torus.  
\endabstract
\endtopmatter
\document

\smallpagebreak

{\bf 0. Introduction.} Throughout $S$ will denote a closed, connected, oriented surface of genus $g \ge 2$. The {\it curve complex} of $S$, denoted $C(S)$, will be the complex whose vertices are the isotopy classes of essential simple closed curves in $S$, and where distinct vertices $x_0,x_1, \dots ,x_k$ determine a $k$-simplex of $C(S)$ if they are represented by pairwise disjoint simple closed curves. If we fix a hyperbolic metric on $S$, then each isotopy class contains a unique geodesic. Moreover two isotopy classes have disjoint representatives if and only if their geodesic representatives are disjoint. We will thus always think of vertices as being geodesics and will use the same notation for a simplex of $C(S)$, the corresponding collection of mutually exclusive simple
closed curves in $S$, and their union as a subset of $S$. 

A simplex $X$ of $C(S)$ determines a {\it compression body}
$$
V_X = S \times [0,1] \cup_{X \times 1} 2-\text{handles } \cup 3-\text{handles}
$$
obtained by attaching $2$-handles along the components of $X \times 1$ and filling in any resulting $2$-sphere boundary components with $3$-cells. $S \times 0$ is called the {\it outer boundary} of $V_X$ and is naturally  identified with $S$.

A pair $X,Y$ of simplexes of $C(S)$ determine a (Heegaard) {\it splitting}
$$
(S;V_X,V_Y) \text{ of a 3-manifold } M_{X,Y} =V_X \cup_SV_Y
$$
for which $(S; X,Y)$ is a (Heegaard) {\it diagram}.

Our goal is to study (compact, oriented) 3-manifolds and  their splittings in terms of the geometry and combinatorics of $C(S)$. We will be primarily interested in the case of closed 3-manifolds ($V_X$ and $V_Y$ handlebodies).

There is a subcomplex $K_X \subset C(S)$ consisting  of those  simplexes $X'$ (and their faces) with $(V_{X'}, S) = (V_X,S)$ (see Lemma 1.2). So the pair $K_X, K_Y$ of subcomplexes of $C(S)$ describes the different diagrams for a fixed splitting. The major questions: What can one say about a splitting in terms of a representative diagram? What can one say about a 3-manifold in terms of the complexes associated with a given splitting?

Perhaps the most natural thing to consider is the {\it geodesic
distance function} $d$, defined on the $0$-skeleton of $C(S)$ by
$d(x,y) = $ the minimal number of $1$-simplexes in a simplicial path
joining $x$ to $y$. So $d(x,y) \le 1$ if and only if $x \cap y =
\emptyset$ and $d(x,y) \le 2$ if and only if there is some $z$ with $x
\cap z = y \cap z = \emptyset$ (i.e $x \cup y$ does not {\it fill} $S$).  Higher distances are harder to visualize, but it is known [MM] that $C(S)$ has infinite diameter with respect to $d$. The significance of the distance function to 3-manifolds begins with:

\proclaim{Observation} The splitting $(S;V_X,V_Y)$ of a closed, oriented 3-manifold $M_{X,Y}$ is:

(1) reducible if and only if $d(K_X,K_Y)=0$, and

(2) weakly reducible if and only if $d(K_X,K_Y) \le 1$.
\endproclaim

Here the distance $d(K_X,K_Y)$ is the minimal distance between their respective vertices; we call it the {\it distance of the splitting}.  The above observation is merely a restatement of definitions in terms of the distance on $C(S)$. Its significance lies in the theorems of Haken [H] that a  splitting without any  cancelling handle pairs is reducible if and only if the corresponding  manifold contains an essential 2-sphere and of Casson and Gordon [CG1] that a weakly reducible splitting is either reducible or the corresponding manifold contains an incompressible surface. Since splittings of $S^3$ are standard [W], a reducible splitting of an irreducible 3-manifold must have a cancelling pair of handles.

We show  in section 3:

\proclaim{Theorem} Let $M$ be a closed oriented 3-manifold which is Seifert fibered or which contains an essential torus. Then any splitting of $M$ is a distance $\le 2 $ splitting.

\endproclaim

The converse is false. As we observe in section 1 there are many hyperbolic 3-manifolds with  distance 2 splittings. For example, any Dehn surgery on a 2-bridge knot (most of which are  hyperbolic manifolds) has a distance $\le 2$ splitting. This follows from

\proclaim{Theorem} If $M$ is obtained by surgery on a link $L$ in $S^3$ then any splitting of $M$ which is derived from a bridge presentation of $L$ is a distance $\le 2$ splitting.
\endproclaim

Of course, these splittings need not be irreducible nor of minimal genus.

However it is true that:

\proclaim{Theorem} There  are distance $n$ splittings of closed,
oriented 3-manifolds for arbitrarily large $n$.
\endproclaim

This is shown in section 2 with an argument supplied by Feng Luo. We
also give in section 5 an explicit construction of some distance $\ge
3$ splittings. However we are unable to answer:

\proclaim{Question} For each $n \ge 3$ are there closed, oriented 3-manifolds which have no irreducible ( or no  minimal genus) splittings of distance $< n$ ?
\endproclaim

We remark that there are 3-manifolds with inequivalent minimal genus splittings [M]  and 3-manifolds with irreducible splittings of different genus [CG2], but any two splittings of a given 3-manifold are stably equivalent [R]. However adding a cancelling pair of handles reduces the distance of the splitting to  zero. It is not clear whether distance survives to any sort of meaningful invariant for 3-manifolds.

In section 2 we introduce some estimates on the distance function which allows us to prove

\proclaim{Theorem} $diam(K_X) = \infty$.
\endproclaim

Which gives an independent proof that $diam(C(S))=\infty$.

This also indicates why the problem is difficult -- one can have ``simple'' splittings represented by diagrams $(S;X,Y)$ with $d(X,Y)$ arbitrarily large. However, we show in sections 3 and 4 that there are obstructions, computable from a fixed diagram, for the corresponding splitting and/or manifold to be reducible, weakly reducible, Seifert fibered, contain an essential torus, or  be a distance 2 splitting.
Examples of their application are given. They arise from enumerating the ``square'' regions of $S - X \cup Y$ according to where the edges lie, and are encoded in a {\it stack intersection matrix}. This turns out to be a much more accurate measure of the real complexity of the splitting. This builds on ideas introduced by Casson and   Gordon [CG2] and extended by Kobayashi [K]  as an obstruction to being weakly reducible. They also provide lower bound estimates for some natural invariants of splittings such as the minimal intersection number between essential disks in the two halfs of the splitting.

Section 6 gives an analysis of all genus two, distance two
splittings. I have become aware of a manuscript by A. Thompson [T]
which analyzes such splittings and includes a proof of Corollary 3.7
(torus case).

\medpagebreak

{\bf 1. Preliminaries.} Throughout $S$ will denote a closed, connected, oriented surface of genus $g \ge 2$. The {\it geometric intersection number} of simple closed curves $J_1, J_2$ in $S$ is
$$
i(J_1,J_2) = min\{\#(J'_1 \cap J'_2): J'_i \text{ isotopic to } J_i \}.
$$
 A collection ${J_1, J_2, \dots }$ of simple closed curves will be said to meet {\it efficiently} if they are in general position and $i(J_i, J_j)=\#(J_i \cap J_j)$ for all $i,j$. This is equivalent to having no disk $D \subset S$ with $\partial(D) = a_i \cup a_j$ where $a_i$ and $a_j$ are arcs in $J_i$ and $J_j$
 respectively. 

If $N$ is a codimension one, bicollared submanifold of a manifold $M$, the result of {\it splitting } $M$ along $N$ will be a manifold $M^*$ whose boundary contains disjoint copies $ N^+$ and $ N^-$ of $N$ together with a map $f :M^* \to M$ which maps $M^* - N^+ \cup N^-$  homeomorphically onto $M - N$ and maps each of $N^+$ and $N^-$ homeomorphically onto $N$.

The {\it curve complex} $C(S)$ is the complex whose $k$-simplexes are the isotopy classes of collections of  $k+1$ mutually exclusive, pairwise non isotopic, essential simple closed curves in $S$. $dim(C(S)) = 3g-4$: a principal simplex of $C(S)$ is a collection of $3g-3$ simple closed curves which splits $S$ into {it pairs of pants} (thrice punctured 2-spheres).

 We will not distinguish notationally between simple closed curves and their isotopy classes.

If $X= (x_0,x_1, \dots , x_k)$ is a $k$-simplex of $C(S)$, we define:

$$
N_X = \text{ normal closure of } \{x_0,x_1, \dots , x_k\} \text{ in } \pi_1(S)
$$

and

$$
V_X = S \times [0,1] \cup_{X \times 1} 2-\text{handles } \cup 3-\text{handles }.
$$

Then $V_X$ is a {\it compression body} whose {\it outer boundary}, $S \times 0$,  is naturally identified with $S$ and $N_X = ker\{\pi_1(S) \to \pi_1(V_X)\}$ determines $V_X$ up to homeomorphisms which restrict to the identity on $S$.

A (Heegaard) {\it splitting} of a compact, orientable 3-manifold $M$ is a representation of $M$ as the union of two compression bodies intersecting on the outer boundary of each. So a pair $X,Y$ of simplexes of $C(S)$ determines a splitting $(S;V_X,V_Y)$ of the 3-manifold
$$
M_{X,Y} = V_X \cup_S V_Y.
$$
Every such 3-manifold is represented in this way, but our requirement $g \ge 2$ precludes the standard genus zero and  one representation of $S^3$, Lens spaces, and $S^2 \times S^1$.

We will be  concerned primarily with closed 3-manifolds which will be represented as above with $V_X$ and $V_Y$ handlebodies. So we say $X$ is a {\it full} simplex of $C(S)$ if $V_X$ is a handlebody. The following gives equivalent properties.

\proclaim{1.1 Lemma} For $X = (x_0,x_1, \dots , x_k)$ a simplex of $C(S)$ the following are equivalent:

(i) $X$ is full in $C(S)$,

(ii) Every component of $S - X$ is planar,

(iii) $S - X$ has $k-g + 2$ components, and

(iv) $\pi_1(S)/N_X$ is free of rank $g$.
\endproclaim

\demo{proof} The arguement is standard, but it is helpful to note that for any simplex $X$ of $C(S)$ of dimension $k$ that
$$
\pi_1(S)/N_X = \pi_1(S_{g_1})
* \dots \pi_1(S_{g_c})*F_r
$$
where $S-X$ has $c$ components, $r = k-g+2$, and $\sum g_i = g-r$. \qed
\enddemo

We call the pair $X,Y$ of simplexes a (Heegaard) {\it diagram} for the
splitting \linebreak $(S;V_X,V_Y)$. It has  been more traditional to think of diagrams as being given by the smallest dimension simplex which will determine the corresponding compression body ($= g-1$ for a handlebody), but we find it convenient to allow {\it superfluous} vertices -- those which can be omitted without changing the 
compression body -- as opposed to {\it essential} vertices which cannot be omitted. Specifically, a vertex $x$ of $X$  will be superfluous if there are distinct components of $S-X$ on opposite sides of $x$, at least one of which is planar. In fact we find that our theorems provide the strongest results when applied to maximal dimensional simplexes. See comments 4.8.

The following is an easy consequence of a theorem of Luo [L].

\proclaim{1.2 Lemma} Two $(3g-4)$-simplexes $X,X'$ of $C(S)$ determine the same handlebody, $(V_X, S) = (V_{X'}, S)$, if and only if there is a sequence $X=X_0, X_1, \dots X_n = X'$ of $(3g-4)$-simplexes of $C(S)$ such that $X_{i-1} \cap X_i$ is a full $(3g-5)$-face of each for $i = 1,2, \dots, n$.
\endproclaim

Caution: if we leave off ``full'' the statement becomes true for any two $(3g-4)$-simplexes of $C(S)$ [HT].

 So for $X$ a full simplex of $C(S)$
$$\align
K_X &= \{X': X' \text{ is a face of a full simplex } X'' \text{ with } V_X = V_{X''}\} \\
 &= \{X': N_{X'} \subset N_X\}
\endalign
$$
Is a connected subcomplex of $C(S)$ whose boundary lies in the non full $(3g-5)$-simplexes.

We say that a splitting $(S;V_1,V_2)$ is {\it reducible} (respectively {\it weakly reducible}) if there are essential disks $D_i \subset V_i, \, i=1,2$ with
$\partial(D_1)=\partial(D_2)$ (respectively $\partial(D_1) \cap \partial(D_2) = \emptyset$). A reducible splitting can be written as a connected sum of lower genus splittings or has a canceling pair of handles.  If $X$ and $Y$ are full simplexes and $S-X\cup Y$ is not simply connected then the splitting is obviously reducible.

A {\it wave} relative to $X,Y$ and based in $X$ ($Y$) is an arc in $S$
whose endpoints lie in the same component of $X$ ($Y$), whose interior
misses $X \cup Y$, which lies on the same side of $X$ ($Y$) near its
endpoints,  and which is not parallel to an arc in $X$ ($Y$). A diagram $(S;X,Y)$ will be called {\it generic} if $X$ and $Y$ are full simplexes which meet efficiently, $S-X \cup Y$ is simply connected, and there are no waves relative to $X,Y$.

If some wave $w$, say, based in $X$ lies in a component $P$ of $S$ split along $X$ with at least four boundary components, then we can do surgery along $w$  to replace $X$ by a simplex $X'$ of the same dimension with $V_{X'} = V_X$ and
$i(X',Y) <  i(X,Y)$ -- or else we discover an obvious reduction. There are two choices for the surgery; one will always give a simplification. If $P$ had only three boundary components, this will not work. But then by Lemma 1.1 $X$ would have a proper, full face to which we could apply the above procedure, if appropriate.

If $(S;X,Y)$ is a generic diagram with, say, $dim(X) < 3g-4$, then we can expand $X$ to a simplex $X'$ with $X',Y$ generic and $dim(X') = dim(X)+1$  as follows. Split $S$ along $X$ to get a disjoint union of planar surfaces. Collapse the boundary components to vertices and identify the families of parallel arcs of the split open $Y$ to single edges to get a graph with one component each in a disjoint union of 2-spheres. Some component has at least four vertices. Take a tree in this component which does not contain all the vertices and whose removal does not separate the component. The boundary of a regular neighborhood of this tree, pulled back to S, represents a vertex we can add to $X$ to get a simplex $X'$ with $X', Y$ generic. Call this operation a {\it generic expansion}. Together these observations prove:

\proclaim{1.3 Lemma} A pair of full simplexes of $C(S)$ can either  be modified to a generic pair of $(3g-4)$-simplexes which determine the same splitting by

(i) taking full faces,

(ii) surgery along a wave, and

(iii) generic expansion

or one will discover an obvious reduction to the associated splitting in the process.
\endproclaim

Suppose $X$ and $Y$ are faces, not necessarily full, of simplexes $X'$ and $Y'$ of $C(S)$. Then there is a natural inclusion $M_{X,Y} \subset M_{X',Y'}$. If the components of $\partial(M_{X,Y})$ are all tori, then the components of $M_{X',Y'} - M_{X,Y}$ are solid tori and $M_{X',Y'}$ is obtained by Dehn filling on $M_{X,Y}$. Moreover every Dehn filling of $M_{X,Y}$ is obtained in this way: the meridians of the filling solid tori can be isotoped to miss the 2-handles and then pulled back to curves in $S$ to represent additional vertices for $X'$ and/or $Y'$. Clearly $d(K_{X'},K_{Y'}) \le d(K_X,K_Y) \le d(X,Y)$.

Let $L$ be a link in $S^3$. We can always isotope $L$ so that for some 3-ball $B \subset S^3$ $L \cap Int(B)$ is the disjoint union of arcs $b_1,b_2, \dots, b_n$ which cobound mutually disjoint disks $D_1,D_2,\dots,D_n$ in $B$ with arcs in $\partial(N)$, and $L -\cup b_i$ is the disjoint union of arcs $a_1,a_2, \dots , a_n$ in $\partial(B)$. This is called an {\it n-bridge presentation} of $L$, and the  minimal such $n$ is called the {\it bridge number} of $L$.

We can choose a regular neighborhood $N = N(L)$ so that $V = Cl(B-N)$
is a genus $n$ handlebody and $N \cap \partial(B)$ is the disjoint
union of $n$ disks $E_1,E_2, \dots, E_n$, each containing some
$a_i$. Then $Cl(S^3-N)$ is homeomorphic to the result of adding
2-handles to $V$ along any $n-1$ of the curves $\partial(E_1),
\partial(E_2), \dots, \partial(E_n)$. See Figure 1.

\bigpagebreak
\epsfysize=3.16in
\centerline{\epsffile{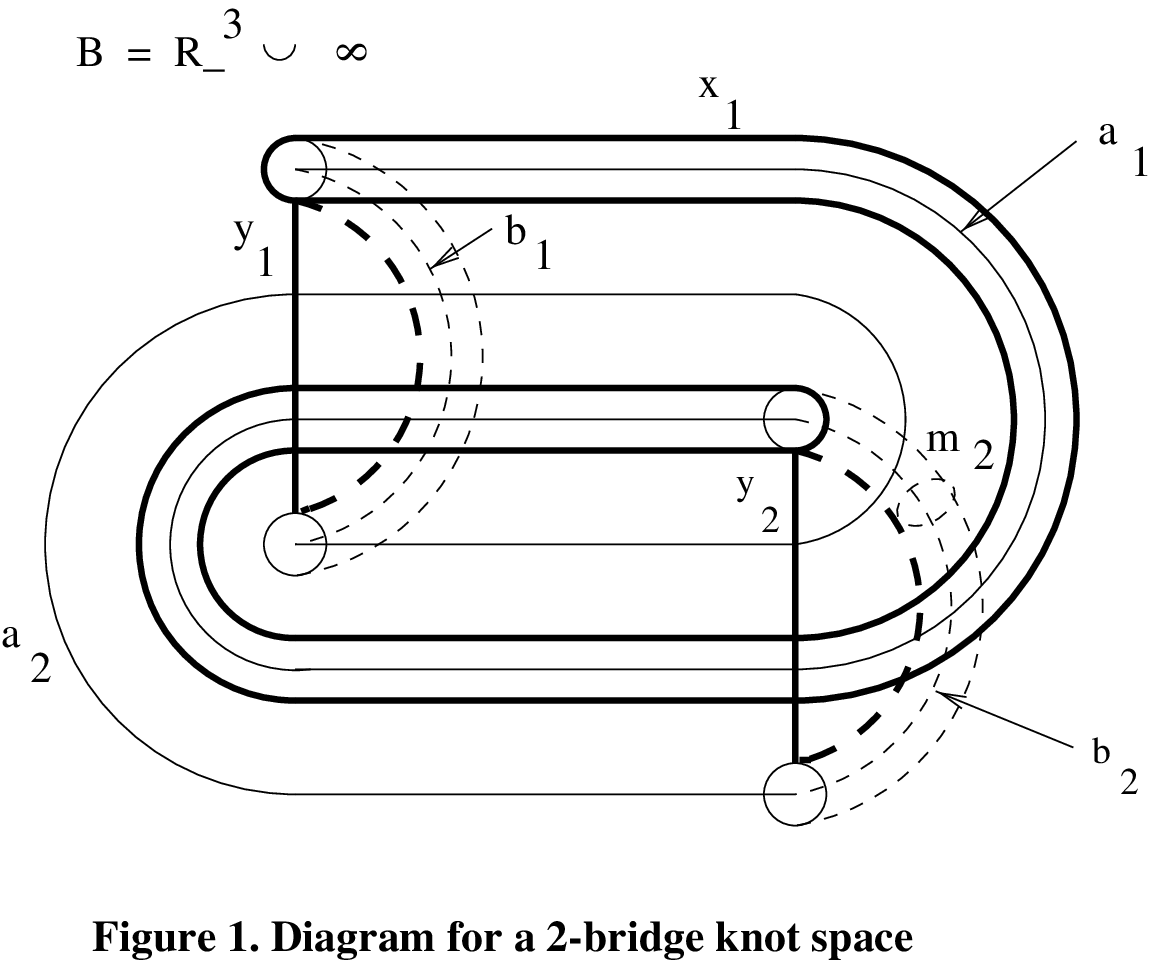}}
\bigpagebreak

So $(S;X,Y)$ is a genus $n$ diagram for a splitting of $Cl(S^3-N)$ where $S = \partial(V)$, $X=\{\partial(E_1), \partial(E_2), \dots, \partial(E_{n-1}\}$, and
$Y= \{\partial(D_1 \cap V), \partial(D_2 \cap V), \dots, \partial(D_n \cap V)\}$. If $M$ is obtained by Dehn surgery on $L$ then this diagram/splitting extends to a diagram/splitting for $M$ which we say is {\it derived from the bridge presentation} of $L$.

\proclaim{1.4 Theorem} If the closed 3-manifold $M$ is obtained by surgery on a nontrivial link $L$ in $S^3$ then any splitting of $M$ derived from a bridge presentation of $L$ is a distance $\le 2$ splitting.
\endproclaim
\demo{Proof} Take a diagram $(S;X,Y)$ for the link complement as described above. So $d(X,Y)$ is an upper bound for the distance of the derived splitting of $M$.

But $d(x_1,y_1) \le 2$; since a meridian $m_2$ (for $S^3-V$) dual to $y_2$ can be chosen disjoint from $x_1 \cup y_1$. \qed
\enddemo

\proclaim{1.5 Corollary} Each irreducible 3-manifold obtained by surgery on a 2-bridge knot has a distance two  genus two  splitting.
\endproclaim

\demo{Proof} A weakly reducible, genus two splitting is easily seen to be reducible; so any irreducible splitting derived from a 2-bridge presentation of the knot will be a distance two splitting. \qed
\enddemo

\medpagebreak

{\bf 2. Distance estimates.} In this section we give upper and lower estimates on the distance between two curves in $C(S)$ which will be used in later sections as well as in establishing that $K_X$ has infinite diameter.

One can easily construct curves at distance two in $C(S)$ which
intersect as much as desired. However the intersection number $i(x,y)$
does provide an upper bound to the distance between curves $x$ and
$y$. To see this note that if one replaces the arc on x between two
points of $x \cap y$ which are adjacent on $y$ by this arc on $y$ one
gets a curve $x'$ which meets $x$ at most once (and so $d(x,x') \le
2$) and which (for appropriate choice) meets $y$ at most half as much as $x$ does. This provides the basis for an inductive proof of:

\proclaim {2.1 Lemma} For vertices $x,y $  of $C(S)$ with $i(x,y) > 0$

$$
d(x,y) \le 2 +  2 log_2(i(x,y)).
$$

\endproclaim

The opposite bound is based on the observation that intersections between curves  which persist on passage to covering spaces have a greater infuuence on their distance. To this end we say that a covering space $p:\tilde S \to S$
{\it separates } simple closed curves $x$ and $y$ in $S$ if there are components $\tilde x$ of $p^{-1}(x)$ and $\tilde y$ of $p^{-1}(y)$ with $\tilde x \cap \tilde y = \emptyset$. A finite covering $p: \tilde S \to S$ is called {\it sub-solvable} if $p$ can be factored as a composition of cyclic coverings (regular with cyclic covering group: which may be assumed to have prime degree).

\proclaim{2.2 Definition} For distinct vertices $x, y$ of $C(S)$ we define the 
{\it covering distance} between $x$ and $y$ to be:

$cd(x,y) = 1+ $ min$\{n:$there is a degree $2^n$ sub-solvable covering of $S$ which separates $x$ and $y \}$.
\endproclaim

\proclaim{2.3 Lemma} Let $x$ and $y$ be distinct vertices of $C(S)$. Then

(i) $d(x,y)=2$ if and only if $cd(x,y)=2$

(ii) $cd(x,y) \le d(x,y)$
\endproclaim

\demo{Proof} Suppose $d(x,y) = 2$ Then $x \cap y \ne \emptyset$ but some vertex $z$ is disjoint from $x \cup y$. We may assume $z$ does not separate $S$; for otherwise $x \cup y$ lies in one component of $S - z$ and we could replace $z$ by a non separating curve in the other component. We construct a double cover of $S$ by glueing together two copies of $S$ split open along $z$. One of these components  contains a (homeomorphic) lift of $x$ and the other a lift of $y$. Thus $cd(x,y)= 2$.

Conversely, suppose $cd(x,y) = 2$ and that $p:\tilde S \to S$ is a double covering separating $x$ and $y$. Then $p^{-1}(x \cup y)$ has two components and these must be interchanged by the non-trivial covering transformation. Some boundary component of a small  regular neighborhood of $p^{-1}(x \cup y)$ projects homeomorphically to an essential simple closed curve in $S -(x \cup y)$.

For (ii) suppose that $d(x,y) = n >2$. Then there are vertices
$$
x=x_0,x_1, \dots , x_n=y
$$
 with $d(x_{i-1},x_i) = 1$. By part (i) there is a double covering
$p^\ast:S^\ast \to S$ separating $x_0 $ and $ x_2$. So for appropriate components
$x^\ast _i$ of $ p^{\ast -1}(x_i)$ we have $x^\ast_0 \cap x^\ast_2 = x^\ast_2 \cap x^\ast_3 = \dots = \emptyset$. So $d(x^\ast_0, x^\ast_n) \le n-1$.

By induction there is sub-solvable covering $q:\tilde S \to S^\ast$  separating $x^\ast_0$ and $x^\ast_n$ and of degree $2^m$ for some $m \le n-2$. Then
 $p= p^\ast \circ q$ has degree $2^{m+1} \le 2^{n-1}$, is sub-solvable, and separates $x$ and $y$. So $cd(x,y) \le m+1 \le n = d(x,y)$. \qed
\enddemo

\proclaim {2.4 Observation} The inequality in (ii) above is, in general,  proper. \linebreak The difference  $d(x,y)-cd(x,y)$ can be made arbitrarily large.
\endproclaim

\demo{Proof} A double cover of $S$ is determined by a homomorphism of $\pi_1(S) $ to ${\Bbb Z}/2{\Bbb Z}$ which in turn is given by the mod 2 intersection number with a fixed curve. The cover in part (i) corresponds to intersection with $z$. Now if (in (ii))  $d(x,y)$ is very large then many of the curves $x_i$ will be homologous (mod 2) to the curve generating the cover, and correspondingly
$d(x^\ast_{i-1},x^\ast_{i+1}) \le 1$ (for appropriate lifts). Thus $d(x^\ast_0,x^\ast_n)$ will be considerably less than $d(x,y)$. However we necessarily have
 $cd(x^\ast_0,x^\ast_n) \ge cd(x,y)-1$. \qed
\enddemo

\proclaim{2.5 Theorem} If $h: S \to S$ is a pseudo-anosov
 homeomorphism and $x$ and $y$ are vertices of $C(S)$ then 
$$
\lim_{n \to \infty}cd(x, h^n(y)) = \infty.
$$
\endproclaim

\demo{Proof} Fix an integer $m >0$ and let $N$ be the intersection of all subgroups of index $2^m$ in $\pi_1(S)$. Let $p:\tilde S \to S$ be the corresponding regular covering space. Then $p$ factors through every degree $2^m$ sub-solvable covering and so any pair of curves which is separated by one of these covers is separated by $p$.

Now $N$ is characteristic and is preserved by $h_\ast$; so $h$ is covered by a pseudo-anosov homeomorphism $\tilde h : \tilde S \to \tilde S$. Now for any essential simple closed curves $z,w$ in $\tilde S$, $lim(i(z,\tilde h^n(w))  = \infty$ -- in fact this property can be taken as a definition of pseudo-anosov (cf [FLP]). Pick a component $\tilde y$ of $p^{-1}(y)$.  Then there is some $n_0$ so that for $n > n_0$,
$\tilde h^n(\tilde y)$ intersects every component of $p^{-1}(x)$. By regularity $p$ cannot separate $x$ and $h^n(y)$ for $n > n_0$. Thus $cd(x,h^n(y)) > m$.
\qed
\enddemo

\proclaim{2.6 Theorem} For $X$ a full simplex of $C(S)$, $diameter(K_X) = \infty$.
\endproclaim

\demo{Proof} One can find two simple closed curves $u, v$ which bound disks in $V_X$, and so represent elements of $N_X$ and such that $u \cup v$ fills $S$.
(cf Figure 2).
The product $h$ of the Dehn twists along $u$ and $v$ is pseudo-anosov [P].
Clearly $h(K_X) = K_X$. The conclusion follows from 2.5 and 2.3. \qed
\enddemo

\bigpagebreak
\epsfysize=2.5in
\centerline{\epsffile{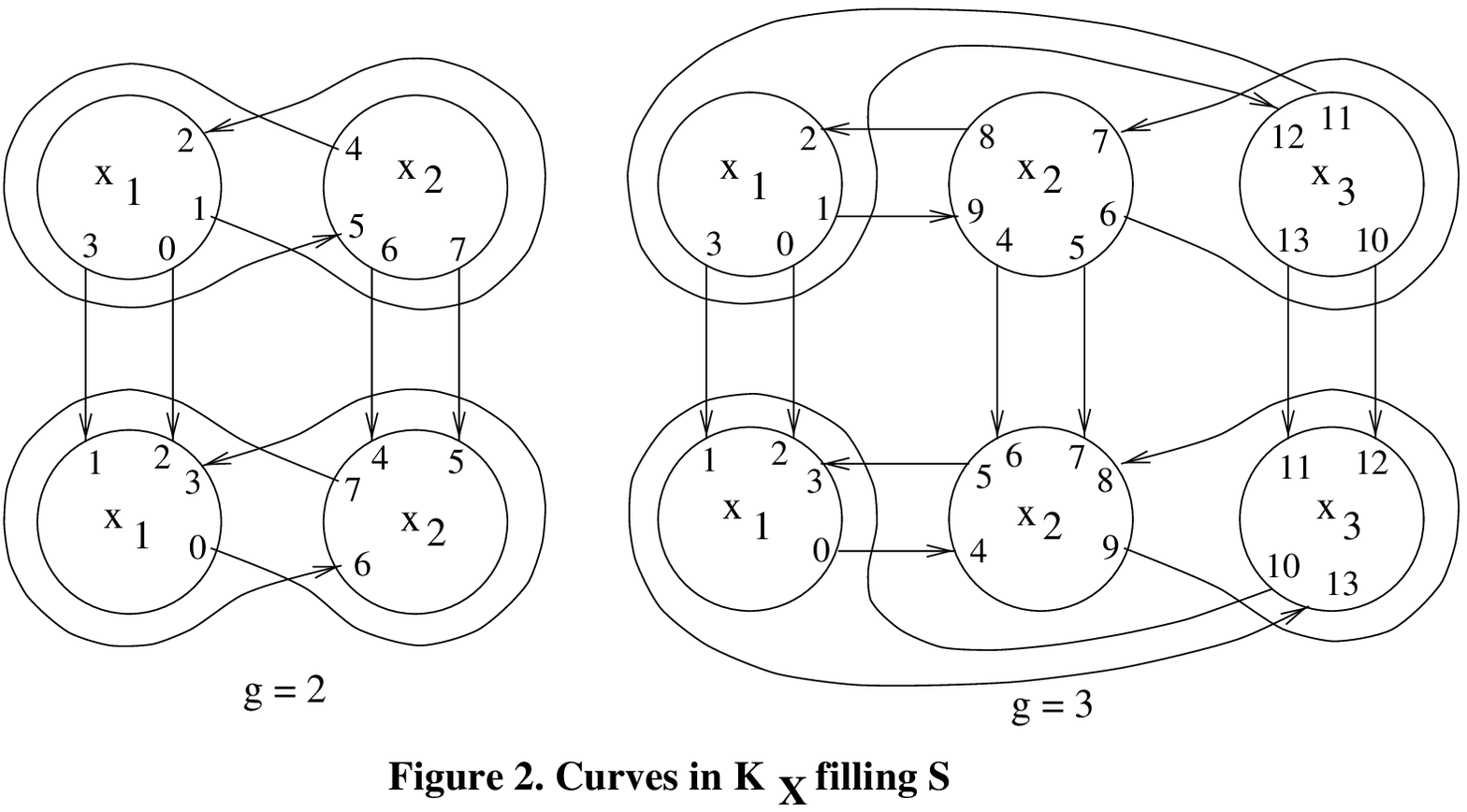}}
\bigpagebreak

A proof for the following was shown to me by Feng Luo who attributes
it to being implicit in the paper [K2] of Kobayashi.

\proclaim{2.7 Theorem} For any $d$ there are full simplexes $X, Y$ of
$C(S)$
with $d(K_X,K_Y) \ge d$.
\endproclaim

\demo{Proof} We regard simple closed curves, with the counting
measure, as elements of the space  $ML(S)$ ($\cong {\Bbb R}^{6g-6}$) of
measured laminations on $S$ and so elements of the space $PL(S)$ 
($\cong S^{6g-7}$) of projective measured laminations on $S$.

Let $X$ be any full simplex of $C(S)$. It is known [Ma] that the
closure, $C$, in $PL(S)$ of the set of vertices of $K_X$ is nowhere
dense in $PL(S)$. So there is a pseudo-anosov homeomorphism $h:S \to
S$ whose stable lamination $L$ is not in $C$. We claim that
$$
\lim_{n \to \infty}d(K_X,h^n(K_X)) = \infty.
$$

To establish this claim it suffices to show that there do not exist
sequences
$$
x_1, x_2, x_3, \dots \text{ and } y_1,y_2, y_3, \dots
$$
of vertices of $K_X$ with $d(x_n,h^n(y_n))$ bounded.

If not then for some $m$ there are sequences
$$
x^j_1,x^j_2,x^j_3, \dots ;j=1,2, \dots, m
$$
of simple closed curves with $x^1_n = x_n, x^m_n = h^n(y_n),$ and $x^j_n
\cap x^{j+1}_n = \emptyset $ for all $n$ and for $j = 1,2, \dots,
m-1$.

Now $i$ extends to a continuous function $i:ML(S) \times ML(S) \to
{\Bbb R}$ and the stable lamination $L$ has non zero intersection with
every lamination which is not a multiple of itself.

By passing to subsequences, we may assume that $x^j_1,x^j_2,x^j_3,
\dots$ converges in $PL(S)$ to some $L_j$ where $L_m = L$. In $ML(S)$
continuity forces $i(L_j,L_{j+1}) = 0; j = 1,2,\dots,m-1$ for any
representatives of the indicated projective classes. By the previous
paragraph and induction we see that $L=L_m=L_{m-1} = \dots =
L_1$. This is impossible since $L_1 \in C$. \qed
\enddemo

\proclaim{2.8 Remark} The above proof also establishes Theorem 2.5
with $d$ in place of $cd$. I do not know the validity of Theorem 2.7
for $cd$.
\endproclaim

\medpagebreak

{\bf 3. Seifert manifolds.} It is known [MS] that every splitting of an orientable Seifert manifold with orientable base space is either {\it horizontal} or {\it vertical} (definitions below). All Seifert manifolds have vertical splittings, but most do not admit horizontal splittings. Theorem 0.3 of [MS] describes, in terms of the Seifert invariants, those  Seifert manifold which have  horizontal splittings.


We show in this section that there is a strong restriction on vertical splittings of Seifert manifolds and a restriction on horizontal splittings as well. An easy corollary is that all splittings of closed, orientable Seifert manifolds are  distance at most two splittings. The same is true for any closed, oriented 3-manifold which contains an essential torus.

\proclaim{3.1  Definition} For a splitting $ \Cal S =(S;V_1,V_2)$ we define the {\it k-simplex intersection complexity} of $\Cal S$ to be:
$$ \split
c_k(\Cal S) & = \text{min}\{i(X_1,X_2): X_i \text{ is a  k-simplex in }Ker(\pi_1(S) \to \pi_1(V_i))\\
&\qquad \qquad \text{ without superfluous vertices}\}
\endsplit
$$
\endproclaim

Perhaps these give the most elementary measures of complexity for a splitting of a 3-manifold. It should be clear adding superfluous vertices to either of the  $X_i$ can only increase intersections (and dimension) without adding topological information, but it should not be assumed that in a fixed dimension
the minimum $i(X_1,X_2)$ occurs without superfluous vertices.

Now consider an orientable Seifert manifold $M$ with base surface $B$ and projection
$f:M \to B$. Suppose we have a cell decomposition of $B$ such that $B
= D \cup E \cup F$ where each of $D,E,F$ is a disjoint union of closed
2-cells of the decomposition, each component of $D$ and of $E$
contains at most one singular point, which is an interior point,  each
component of $F$ is a square containing no singular point and having
one pair of opposite sides in $D$ and the other pair in $E$,  $Int(D)
\cap Int(E) = Int(D) \cap Int(F) = Int(E) \cap Int(F) = \emptyset$,
and $D \cup F$ and $E \cup F$ are connected. See Figure 3.

\bigpagebreak
\epsfysize=1.71in
\centerline{\epsffile{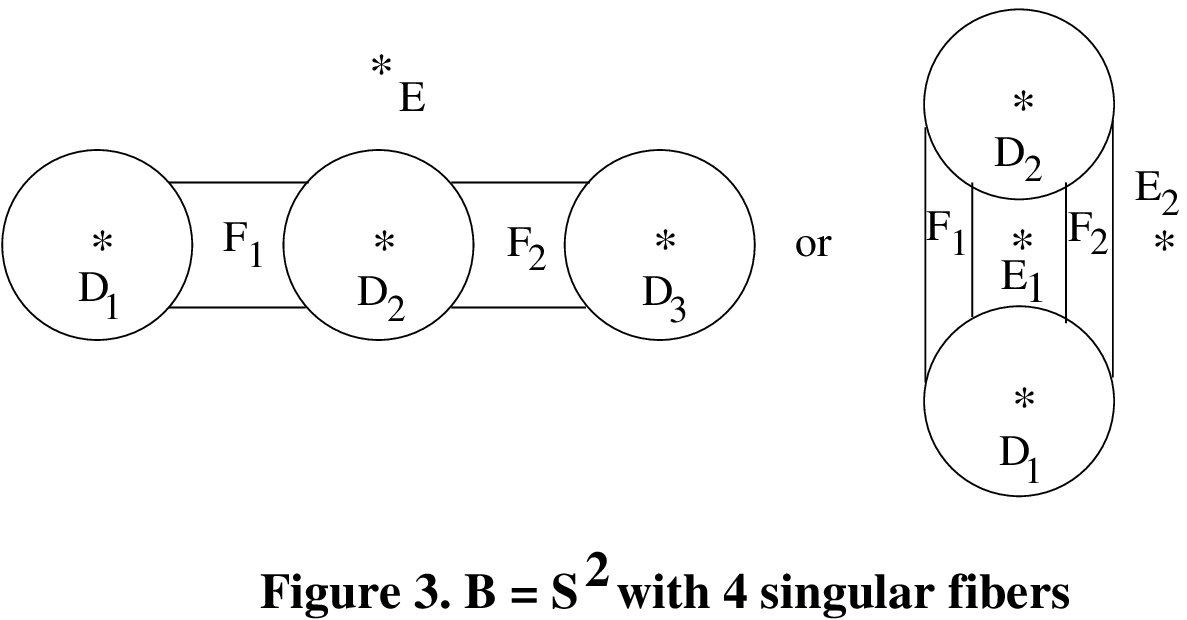}}
\bigpagebreak

Then $f^{-1}(D)$ is homeomorphic to $D \times S^1$,  the same holds
for $E$ and $F$, and we consistently fix such identifications. Let $V_1 = D \times
S^1\cup F \times [0,1/2]$ and $V_2 = E \times S^1 \cup F \times [1/2,1]$; where $S^1 = [0,1]/0 \sim 1$.
Put $S=V_1 \cap V_2 =\partial(V_1) = \partial(V_2)$. Then $(S;V_1,V_2)$ is a splitting of $M$ of genus
$$
g = \beta_0(D)+\beta_1(D \cup F)= \beta_0(E) + \beta_1(E \cup F) = 1 + \beta_0(F)
$$
 which we call a {\it vertical splitting}.

\proclaim{3.2 Theorem} Let $\Cal S $ be an irreducible  splitting of genus $g \ge 2$ of a closed, orientable 3-manifold $M$. Then $\Cal S$ is a vertical splitting of a Seifert manifold if and only if $c_{g-2}(\Cal S) \le 2g-2$.
\endproclaim

\demo{Proof} Given a vertical splitting of $M$ as in the definition above, choose for each component $F_i$ of $F$ spanning arcs $A_i$ and $B_i$ of $F_i$ meeting in a single point and joining the opposite edges of $F_i$ which lie in $D$ and $E$ respectively. Then $X_1= \cup \partial(A_i \times [0,1/2])$ and $X_2=\cup \partial(B_i \times [1/2,1])$ are the desired $g-2$ simplexes with $i(X_1,X_2)=2g-2$.

Now suppose we have a splitting $(S;V_1,V_2)$ of $M$ and for each $i=1,2$ a collection $ D_i =\{D_{i,0}, D_{i,1}, \dots , D_{i,g-2}\}$ of disjoint, properly embedded 2-cells in $V_i$ such that no component of $S-\cup_j \partial(D_{i,j})$ is planar and $i(\partial(D_1),\partial(D_2)) \le 2g-2$.

Suppose $S- \cup_j \partial(D_{i,j})$ has $c_i$ components.
By the reasoning of Lemma 1.1 the sum of their genera is $c_i$. Since no component is planar, they are all of  genus one. Thus $D_i$ splits $V_i$ into a collection of solid tori.

If some $D_{1,j}$ meets some $D_{2,k}$ in a single point, the splitting has a
 trivial handle and is reducible contrary to assumption.

If  some  $D_{1,j}$ does not meet $D_2$  then $\partial(D_{1,j})$ lies in a solid torus component $W$ of $V_2$ split along $D_2$ and we see that $M$ contains a punctured Lens space ( possibly $S^3$ or $S^2 \times S^1$). 
In the case of $S^3$  or $S^2 \times S^1$,  some meridian disk $D$ of $W$ meets $\partial(D_{1,j})$ in at most one point and can be chosen to have boundary in $S$. In the other cases $M$ is either a Lens space or a non trivial connected sum. But Lens spaces have unique splittings [BO] and splittings of connected sums are reducible [H]; so in every case we get the contradiction that the splitting is reducible.

These observations, together with the assumption  $i(\partial(D_1),\partial(D_2)) \le 2g-2$ imply that each $D_{1,j}$ must meet some component of $D_2$, which we assume is $D_{2,j}$ in two points and be disjoint from every other component of $D_2$. Thus a regular neighborhood $N$ of $D_1 \cup D_2$ is the disjoint union of $g-1$ solid tori. The $j$-th one can be identified with $B^2 \times S^1$ in such a way that $D_{1,j}$ is identified with $A_1 \times [0,1/2]$ and $D_{2,j}$ with $A_2 \times [1/2,1]$ for some pair $A_1, A_2$ of properly embedded arcs in $B^2$ which cross once. Each component, say $W$,  of $Cl(V_i -N)$ is a solid torus which  meets $N$ in a collection of annuli in the boundary of each.

The product fibration of $N$ will extend to a Seifert fibration of $W$ unless one of these annuli is inessential in $W$. This, however, would give a punctured Lens space in $M$ and a contradiction as before. Thus we can conclude that $M$ is Seifert fibered and, in fact, see that the given splitting has the structure of a vertical splitting. \qed
\enddemo

\proclaim{3.3 Corollary} A vertical splitting of genus $g$ of a closed Seifert manifold is a distance at most two splitting if $g = 2$ and a distance at most  one splitting otherwise.

\endproclaim

We define a {\it horizontal splitting} as follows. Take a surface bundle
$$
N=F \times [0,1]/(x,0) \sim (\phi(x),1)
$$
where $F \ne B^2$ is a compact, connected, orientable surface with one boundary component and $\phi: F \to F$ is an orientation preserving homeomorphism with $\phi | \partial (F) = 1$. Then $\lambda = \partial(F) \times 0$ and $\mu = x_0 \times [0,1]/ \sim$ form a basis for $H_1(\partial(N))$. Let 
$$
M = N \cup_h B^2 \times S^1
$$ be a Dehn filling of $N$ where $h:\partial(B^2 \times S^1) \to
 \partial(N)$ is a homeomorphism such that $h(\partial(B^2) \times 0)$
 is homologous to
 $\mu + n \lambda$ for some $n \in {\Bbb Z}$. 

Now $h^{-1}(\partial(F)\times \{0,1/2\})$ bounds an annulus $A \subset B^2 \times S^1$ which splits $B^2 \times S^1$ into two solid tori $U_1$ and $ U_2$ with $(U_i, A)$ homeomorphic to $(I \times I \times S^1, I \times 0 \times S^1)$. Then $V_1 = F \times [0,1/2] \cup_h U_1$ and $V_2 = F \times [1/2,1] \cup_h U_2$
are handlebodies which give a splitting of $M$ which we call a {\it horizontal splitting}

If $\phi$ is (isotopic to) a periodic homeomorphism of $F$ then $N$ will be Seifert fibered and if $n \ne 0$ this extends to a Seifert fibration of $M$ (with a singular fiber in $M-N$ when $n \ne \pm 1$). However there is no need to assume this for the following.

\proclaim{3.4 Theorem} A horizontal splitting of a closed 3- manifold is a distance at most two splitting; in fact it has the form $(S;V_X, V_Y)$ where every vertex of $X$ is at distance at most two from some vertex of $Y$ and vice versa.
\endproclaim

\demo{Proof} Let $(S;V_1,V_2)$ be a splitting obtained as in the above
definition. There are homeomorphisms $f_i :F \times [0,1] \to V_i$
such that $f = f_2^{-1} \circ f_1|\partial(F \times [0,1])$ is level preserving, $f(x,0) = (x,0)$, $f(x,1) =( \phi(x),1)$, and $f(e^{i\theta},t) = (e^{in\theta},t)$ for $(e^{i\theta}, t) \in \partial(F) \times [0,1])$.

Let $A$ be a 1-manifold which splits $F$ into a 2-cell. Then $X =
f_1(\partial(A \times [0,1])$ and $Y = f_2(\partial (A\times [0,1]))$
gives the desired diagram for the splitting. A component of $A$ gives
rise to a component of $X$ and a component of $Y$ which are at a
distance of at most two as they don't fill $f_1(F \times 1) = f_2(F
\times 0)$. \qed
\enddemo

\proclaim{3.5 Theorem} Let $(S;V_1,V_2)$ be an irreducible, horizontal splitting of a closed Seifert manifold $M$. Then there are essential annuli $A_i \subset V_i, i = 1,2$ with \linebreak  $i(\partial(A_1), \partial(A_2)) \le 1$
\endproclaim

\demo{Proof} The base surface of $M$, homeomorphic to $F/\phi$, is orientable.
If it has positive genus or if $M$ has more than three singular fibers, then $M$ contains an essential torus which is a union of nonsingular fibers and which may be assumed to lie in $N$ (as in the definition of horizontal splitting). This torus will meet $V_1$ and $V_2$ in essential annuli which may be moved slightly to be disjoint.

 Lens spaces do not have irreducible splittings of genus $g \ge 2$. So we are left with the case that $M$ fibers over $S^2$ with exactly three exceptional fibers.

The horizontal splitting comes from a Dehn filling on a surface bundle \linebreak $F \times I/(x,0) \sim (\phi(x),1)$ for some periodic homeomorphism $\phi:F \to F$. The orbit map $f: M \to S^2$ induces a cyclic branched covering $f:F \to B^2$  with $\phi$ generating the group ( $\cong {\Bbb Z}_n$) of covering transformations, and this, in turn, is determined by an epimorphism $\pi_1(B^2 -\{\text{branch points}\}) \to {\Bbb Z}_n$. 

Now $F$ is obtained from an annulus $A$ by identifying edges, in pairs, on one component of $\partial(A)$ and $\phi$ is induced by an equivariant rotation of $A$. To see this choose an arc $a \subset Int(B^2)$ which contains the branch points. Then $\pi_1(B^2 -a) \to {\Bbb Z}_n$ is an epimorphism (as $\partial(F)$ is connected); so $f^{-1}(B^2 -a)$ is a half open annulus in $F$ which we can complete to get the desired $A$.

Some arc in $A$ whose end points are midpoints of identified edges gives rise to an essential simple closed curve $J$ in $F$. Otherwise we get the contradiction that $F = B^2$. Moreover we can adjust so that $\#(J \cap \phi(J)) = i(J,\phi(J)) \le 1$. Then $A_1 = J \times [0,1/2]$ and $A_2 = J \times [1/2,1]$ (adjusted to general position) give the desired annuli. \qed
\enddemo

\proclaim{3.6 Lemma} Let $(S;V_1,V_2)$ be a strongly  irreducible splitting of a closed 3-manifold $M$ which contains an essential torus or Klein bottle $T$. Then, after an isotopy of $T$ we may assume that each component of $T \cap V_i, i = 1,2$ is an essential annulus or M\"obius band in $V_i$.
\endproclaim

\demo{Proof}. The sum of the Euler characteristics of the components of $T$ split along $S \cap T$ is $\chi(T) = 0$. So if there are no disk components, then all components are annuli or M\"obius bands which are incompressible in the $V_i$ in which they lie. Some of these annuli might be parallel to annuli in $S$ and  could be eliminated by an isotopy of $T$, but some must remain; as $T$  cannot be isotoped into a handlebody.

So we induct on the number of disk components of $T$ split along $S \cap T$. Since the splitting is assumed to be strongly irreducible, there are not disk components in both $V_1$ and $V_2$. We assume they are all in $V_1$. Choose a component $C$ of $T\cap V_2$ such that $\partial (C)$ contains a simple closed curve $J$ which bounds a disk $D \subset T \cap V_1$. $J$ does not bound a disk in $V_2$; so it must meet some meridian disks for $V_2$. This means that there is a boundary compression of $C$ along an arc $a \subset C$ which cobounds a disk in $V_2 -T$ with an arc $b \subset S$ where $a$ has an endpoint in $J$. The boundary compression eliminates $D$ as a component of $T \cap V_1$: replacing it by an annulus, if both endpoints of $a$ are in $J$ or by the connected sum of $D$ with another component of $T \cap V_1$ otherwise.

If $C$ were an annulus it would be reduced to a disk in $V_2$ by the boundary compression. This, however, would contradict strong irreducibility. There is no other way in which new disk components could be introduced; so induction applies to complete the proof. \qed
\enddemo

\proclaim{3.7 Corollary} Any strongly irreducible splitting of a closed 3-manifold $M$ which contains an essential torus or Klein bottle is a distance two splitting.
\endproclaim

\demo{Proof} Let $(S;V_1,V_2)$ be a strongly irreducible splitting of our manifold $M$. By 3.6 there are essential annuli or M\"obius bands $A_i \subset V_i, i = 1,2 $ which have a common boundary component $J$. A boundary compression of $A_i$ gives rise to a disk $D_i \subset V_i$ which, since $A_i$ is not parallel to an annulus in $S$, is essential in $V_i$. We may assume $\partial(D_i) \cap \partial (A_i) = \emptyset$. Thus $d(\partial (D_i), J) \le 1$ and the proof is complete. \qed
\enddemo

\proclaim{3.8 Corollary} Any splitting of a closed, orientable Seifert manifold is a distance at most two splitting.
\endproclaim

\demo{Proof} For Seifert manifolds with orientable base this follows from 3.3, 3.4 and Theorem 0.1 of [MS]. Seifert manifolds with nonorintable base must contain an essential torus or Klein bottle and 3.7 applies. \qed
\enddemo

\medpagebreak

{\bf 4. Complexity bounds.} Casson and Gordon [CG2], [K] gave a {\it rectangle condition} on a Heegaard diagram which implies that it determines a strongly irreducible splitting. We give here a quantitative version of this condition which gives lower bounds for the complexity $c_k(\Cal S)$, and in particular is used, in Theorem 4.4, to show that  a diagram does not determine a vertical splitting of a  Seifert manifold manifold. We also give a quantitative version of the {\it strong rectangle condition} introduced by 
 Kobayashi [K] to give conditions that the splitting determines an atoridal manifold. The quantitative version is used in corollary 4.7 to give criteria that a splitting not be a horizontal splitting of a Seifert manifold. 

Somewhat stronger versions of the  results are available for genus two splittings and are presented separately.

A  pair $X,Y$ of simplexes of  $C(S)$   determines a cell structure on $S$ whose faces are the components of $S -(X \cup Y)$ (assuming these are all simply connected). Every vertex has order four and every face has an even number ($ \ge 4$) of edges which lie alternately in $X$ and $Y$. A standard calculation gives.

\proclaim{4.1 Lemma} If $X$ and $Y$ are simplexes of $C(S)$ with $S-(X \cup Y)$ simply connected and having $n_i$ $2i$-gon components ($i = 1,2, \dots$), then
$$
\chi(S) = \sum(1-i/2)n_i
$$
\endproclaim

Since $n_1 = 0$ (assuming efficient intersection) and $\chi(S) < 0$, most of the complementary regions will be squares with one pair of opposite edges in $X$ and the other pair in $Y$.
If we ``stack'' together adjacent squares along common edges in $X$ maximally we get an {\it $X$-stack}. The {\it top} and {\it bottom} edges will lie in {\it large} ( $ \ge 6$ sides) regions, and the {\it sides} will lie in (possibly the same) component(s) of $Y$. The process must actually stop at a top and bottom; otherwise we would have two parallel components of $Y$. The {\it $Y$-stacks} are defined by interchanging the roles of $X$ and $Y$. The number of squares in a stack is called its {\it height}. For logical consistancy we must include stacks of height 0 -- corresponding to edges common to two large complementary regions. These occur rarely and will never satisfy the  conditions of our theorems.

For a somewhat different picture, split $S$ open along $Y$. $X$ gets split into a collection of arcs which fall into  familys of parallel arcs which correspond to the $X$-stacks; where a family with $h+1$ arcs corresponds to a stack of height $h$. A component of $Y$ will contain sides of some $X$-stacks which lie to either of its sides; an $X$-stack on one side may meet several $X$-stacks from the other side. The union of the $X$-stacks is a regular neighborhood of a train track which has one branch for each $X$-stack. This train track carries $X$ with weights the numbers: stack height +1.

If $X$ and $Y$ are full simplexes forming a generic pair. Then each $X$-stack lies  in a component $P$  of $S$ split along $Y$, and $P$ is planar. The sides of the stack must lie in different boundary components of $P$ (but which could be identified to the same, essential,  component of $Y$); otherwise  there would be a wave. In fact every potential wave (arc in $P$ with ends in the same boundary component of $P$ which is not parallel to an arc in this boundary component) must cross some $X$-stack (possibly of 0 height).
This gives part of: 

\proclaim{4.2 Lemma} If $X,Y$ is a generic pair of full simplexes of $C(S)$ then there are the same number of $X$-stacks as $Y$-stacks. This number is $\sum n_i/2$ and it lies in the interval  $[2max\{dim(X),dim(Y)\}+2, 6g(S)-6]$.
\endproclaim

\demo{Proof} The function that assigns to each $X$-stack its bottom (with respect to some arbitrary orientation) gives a bijection between the set of $X$-stacks and half of the edges of large regions. This is clearly symmetric in $X$ and $Y$.

When we split $S$ along $Y$ we get $dim(Y)+2-g$ component, each planar, with a total of $2dim(Y)+2$ boundary components . It takes at least $p$ stacks to block a wave in a planar region with $p$ boundary components.This gives the lower bound

If we collapse the boundary components of $S$ split along $Y$ to points and collapse the $X$-stacks to arcs, we get a cell structure on the disjoint union of $dim(Y)+2-g$ $2$-spheres ($g=g(S)$)   whose edges correspond to the $X$-stacks and whose number is at most the number of edges  in a triangulation with the same number of vertices, which is $6g-6$. \qed
\enddemo

The {\it intersection number} of an $X$-stack and a $Y$-stack is the number of squares common to the two stacks. The {\it stack intersection matrix} for a pair $X, Y$ is the matrix of intersection numbers of non empty (but possibly height 0) $X$-stacks and $Y$-stacks.
A  stack of height 0 will give a row or column of 0's to this matrix. The matrix will have $\sum n_i/2$ rows and columns. If all the large complementary regions are hexagons, this number will be $6g-6$. This holds, in particular,  when $X$ and $Y$ are maximal dimension ($=3g-4$) simplexes forming a generic pair. Then  each component of $S$ split along $Y$ contains three $X$-stacks: one with sides in each pair of boundary components (and vice versa).

 
The stack intersection matrix defines an integer valued bilinear form which we can use to estimate intersection numbers as follows. Suppose $A$ is a simplex of $C(S)$ which has been isotoped so as to meet $X$ and $Y$ efficiently. If $t_j$ is an $Y$-stack, a component of $A \cap t_j$ which does not meet the sides of $t_j$ must have one end point in the top and one in the bottom of $t_j$. We call such a component a {\it stack crossing} and denote the number of such by $a_j \ge 0$. So we get a {\it  stack crossing vector} $(a_1, a_2, \dots )$ of $A$ with respect to $Y$. If another simplex $B$ of $C(S)$ has stack crossing vector $(b_1,b_2, \dots )$ with respect to $X$ then in each square common to $X$- stack $s_i$ and $Y$-stack $t_j$ we see $b_ia_j$ points of $A \cap B$. This gives

\proclaim{4.3 Lemma} Let $X,Y$ be a generic pair of simplexes of $C(S)$ with corresponding stack intersection matrix $(u_{i,j})$.  Suppose $A$ and $B$ are simplexes of $C(S)$ which meet 
efficiently and which have stack crossing vectors $(a_1,a_2, \dots )$
with respect to $Y$ and $(b_1,b_2, \dots )$ with respect to $X$. Then
$$
i(A,B) \ge \sum_{i,j} b_iu_{i,j}a_j
$$
\endproclaim

\proclaim{4.4 Theorem} Let $X,Y$ be a generic pair of simplexes of $C(S)$ with corresponding stack intersection matrix $(u_{i,j})$. Then $c_0(S;V_X,V_Y) \ge
4 min\{u_{i,j}\}$; so if $u_{i,j} > 0$ for all $i,j$ the splitting $(S;V_X,V_Y)$ is strongly irreducible and is not a vertical splitting of a Seifert manifold.
\endproclaim
\demo{Proof} Let $a$ and $b$ be vertices of $K_X$ and $K_Y$ respectively. We may assume that all pairs of curves from $X,Y,a,b$ meet efficiently. Suppose that $a \cap X \ne \emptyset$. Look at how  a disk in $V_X$ bounded by $a$ meets disks bounded by the components of $X$. There must be at least two ``outermost'' arcs of intersection. The arcs on $a$ with the same boundary must contain $Y$-stack crossings -- otherwise there would be a wave. If $a \cap X = \emptyset$, then $a \cap Y \ne \emptyset$ -- otherwise $S - X \cup Y$ is not simply connected. Thus $a$ meets some $Y$-stack in a stack crossing. If there were only one $Y$- stack crossing on $a$ we could replace a neighborhood of this crossing on $a$ by a pair of arcs running parallel to the top and bottom of the stack and ending in a component of $X$. There are two ways of doing this. One will produce a wave.

Thus $a$ contains two $Y$-stack crossings (possibly of the same stack. Similarily $b$ contains two $X$-stack crossings. It follows from Lemma 4.3 that $i(a,b) \ge 4 min\{u_{i,j}\}$. The final conclusion then follows from Theorem 3.2 \qed
\enddemo

This theorem can be improved for genus two splittings due to the fact that $a$ and $b$ must then contain crossings of two distinct stacks (though not necessarily associated with outermost arcs). So let $g(S)=2$ and let $X=(x_0,x_1,x_2)$ and $Y=(y_0,y_1,y_2)$ be
 generic $2$-simplexes in $C(S)$ without essential vertices. 

Each $x_i$ ($y_j$) misses two $Y$-stacks ($X$-stacks): those with sides in the other two components of $X$ ($Y$). Let $c_{i,j}$ be the sum of the corresponding four stack intersection numbers.

\proclaim{4.5 Theorem} Let $X$ and $Y$ be a generic pair of $2$-simplexes in $C(S)$ without essential vertices; where $g(S) =2$. Then $c_0(S;V_X,V_Y) \ge min\{c_{i,j}\}$.
\endproclaim

\demo{Proof} Choose vertices $a \in K_X, b \in K_Y$ as before. Consider first the case $a \cap X \ne \emptyset \ne b \cap Y$. Now $a$ meets each of the (pants) components $P_1, P_2$ of $S$ split along $X$ in at most three families of parallel arcs. Corresponding to an outermost arc, one of these, say $P_1$, contains a family of $n \ge 1$ parallel arcs in $a$ whose end points lie in the same component, say $x_i$ of $X$. Consider the various possibilities for the families of arcs in $a \cap P_2$ and note that there are matching equations equating the number of endpoints on each component of $X$ comming from the two sides. There is only one solution (this is the basis of Neilsen-Fenchel coordinates, c.f [P2]): $a$ must meet $P_2$ in a family of $n$ parallel arcs with both end points in $x_i$.

Similarly we see that for some $j$ and some integer $m \ge 1$ that $b$ must meet each component of $S$ split along $Y$ in families of $m$ parallel arcs with both end points in $y_j$. It then follows from Lemma 4.3 that $i(a,b) \ge nmc_{i,j}$.

Now consider the case $a \cap X = \emptyset$ but with $b$ as above. Then $a$ is isotopic to a component of $X$, say $x_0$. We can isotope $a$ to either side of $x_0$ and each time  get a lower bound on $i(a,b)$ from Lemma 4.3. Averageing these gives $i(a,b) \ge m(c_{1,j} + c_{2,j})/2$.

Simerlarly one can show, for example, that $i(x_0,y_1) \ge (c_{1,0}+c_{1,2}+c_{2,0}+c_{2,2})/4 $ to handle the case $a \cap X = b \cap Y = \emptyset$. \qed
\enddemo

We will say that a pair $X, Y$ of simplexes of $C(S)$ is {\it
complete} if they are full simplexes,  it  is a generic pair, for each
pair $x_i, x_j$ of components of $X$ which lie in the closure of a
component $U$ of $S-X$ there is a $Y$-stack in  $\bar U$ with one side in
$x_i$ and the other in $x_j$, and the symmetric condition holds as well. We
require that this holds for $x_i = x_j$ when $U$ contains both sides
of this curve

\proclaim{4.6 Theorem} Let $X, Y$ be a  pair  of $(3g-4)$-simplexes of
$C(S)$ ($g = g(S)$) such that  each pair $X',Y'$ of  $3g-5$ faces $X'$
of $X$ and $Y'$ of $Y$ is complete and has all stack intersection numbers at least $u$. 
If $A$ is an essential annulus in $V_X$ and $B$ is an essential annulus in $V_Y$, then $i(\partial (A), \partial (B) \ge u$.
\endproclaim

\demo{Proof} The goal is to find $3g-5$ faces $X'$ of $X$, $Y'$ of $Y$ such that $\partial (A)$ contains an $Y'$-stack crossing and $\partial (B)$ contains a $X'$-stack crossing and then apply Lemma 4.3. A crossing of a $Y$-stack by $\partial (A)$ (or a $Y$-stack by $\partial (B)$ will give one for every codimension one face ; so we assume none exist.

The intersection of $A$ with the union $D$ of the disks bounded by the components of $X$ will be a collection of arcs splitting $A$ into 2-cells. The graph in $A$ dual to these arcs is a deformation retract of $A$ and thus has Euler characteristic zero. If this graph has a vertex of order different than $2$, it must have a vertex of order $1$. This vertex comes from an outermost arc and indicates the presence of an $Y$-stack crossing in $\partial (A)$.

So we assume all the vertices of the dual graph have order $2$ and thus the components of $A$ split along $D$ will all be squares with a pair of opposite sides in $D$ and the other pair in $S$. The first pair must lie in different components of $D$; otherwise there would be an $Y$-stack crossing. Each such square lies in some component $C$ of $V_X$ split along $D$. $C$ is a 3-cell which meets $D$ in three 2-cells in its boundary. The square misses exactly one of these 2-cells and this gives a preferred direction in which to isotope the square into $S$ -- away from the 2-cell it misses. If these directions agree from square to square, we get the contradiction that $A$ is parallel to an annulus in $S$.

Thus there must be adjacent squares $s_1, s_2$ in adjacent components $C_1, C_2$ of $V_X$ split along $D$ such that $\partial (s_1 \cup s_2)$  separates some pair $x_i, x_j$ of components of $X$ on $\partial (C_1 \cup C_2)$ (see Figure 4).

\bigpagebreak
\epsfysize=1.75in
\centerline{\epsffile{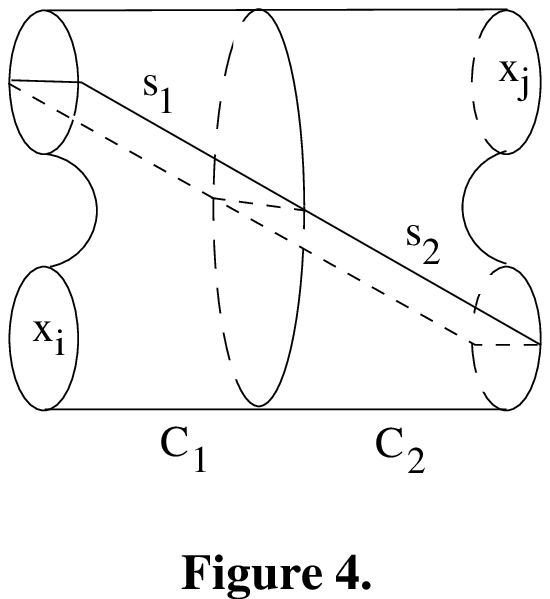}}
\bigpagebreak

 We get the desired face $X'$ by deleting the component of $X$  lying in $C_1 \cap C_2$. For any codimension one face $Y'$ of $Y$ there is, by assumption, a $Y'$-stack with sides in $x_i$ and $x_j$ which therefore must be crossed by $\partial (A)$. \qed
\enddemo

\proclaim{4.7 Corollary} With the hypothesis and notation of Theorem 4.6, if $u > 0$ then $M_{X,Y}$ contains no essential torus. If $u > 1$ then
$(S;V_X,V_Y)$ is not a horizontal splitting of a Seifert manifold.
\endproclaim

\demo{Proof} This follows immediately using  3.5 and 3.6 and the observation that the splitting is strongly irreducible by 4.4. \qed
\enddemo

\proclaim{4.8 Comments}

1. Suppose $(S:X',Y)$ is obtained from $(S;X,Y)$ by a generic expansion. The square regions of $S-X' \cup Y$ will be ``halves'' of squares of $S-X \cup Y$ split by the new curve together with some new squares cut off of large regions of $S -X \cup Y$ by the new curve. A new square could form a $X'$-stack of height one (if two large regions of $S- X \cup Y$ share an edge), but in general they get added to (subdivided) $X$-stacks and one can show that the stack intersection numbers go up. For this reason the theorems of this section generally give the best results when applied to principal simplexes $X,Y$ -- that is pants decompositions.

2. Theorem 4.5 is definitely stronger than Theorem 4.4 applied to genus two splittings. Consider the example of Figure 10(a). Even when extended to pants decompositions, some stack intersection numbers will be zero, but the minimal $c_{i,j}$ will be $2$.
\endproclaim

\medpagebreak

{\bf 5. Distance three splittings.} In this section we give a criterion for recognizing distance three splittings and apply it to give some examples.

For simplexes $X$ and $Y$ of $C(S)$, $s$ an $X$-stack and $t$ a $Y$-stack, $s \cup t$ and \linebreak $Cl(S - s \cup t)$ will be  2-manifolds except possibly at a finite number of singular points where a corner of $s$ meets a corner of $t$. A regular neighborhood of $Cl(S - s \cup t)$ will be called a {\it complementary region} and denoted $CR(s \cup t)$. One gets one  such by adding to $Cl(S - s \cup t)$ a suitable neighborhood of the singular points. We say that $s \cup t$ {fills} ({\it almost fills}) $S$ if the components of $CR(s \cup t)$ are all 2-cells ( 2-cells and annuli)

\proclaim{5.1 Theorem} Let $X, Y$ be a generic pair of full simplexes of $C(S)$
without essential vertices such that for every $X$-stack $s$ and $Y$-stack $t$ $s \cup t$ fills $S$. Then $(S;V_X,V_Y)$ is a distance  $ \ge 3$  splitting.
\endproclaim

\demo{Proof} Suppose there are vertices $a \in K_X, b \in K_Y$ with $d(a,b) \le 2$. So there is some vertex $c \in C(S)$ such that $c \subset S - a \cup b$. Now $a$ must cross some $Y$-stack $t$ and $b$ must cross some $X$-stack $s$. Then $c$ misses the stack crossings and therefore must be  isotopic into $CR(s \cup t)$.  This fact requires the assumption
 that $X$ and $Y$ have no essential vertices; so that the stacks are embedded rectangles. In particular, if two squares of $s \cup t$ intersect in an edge lying in $Y$ then at most one of the squares can lie in $X$ and vice versa. It follows from the observation that the regions of $s \cup t$ complementary to the crossing arcs either lie in the interior of one of the stacks or is a rectangle or annulus (it has a cell decomposition by rectangles) meeting $\partial(s \cup t)$ in a connected set. It of course contradicts the assumption that $s\cup t$ fills $S$. \qed
\enddemo

\proclaim{5.2 Theorem} Let $X, Y$ be a generic pair of full simplexes of $C(S)$
without essential vertices such that for every $X$-stack $s$ and $Y$-stack $t$ $s \cup t$ almost fills $S$. Then either $(S;V_X,V_Y)$ is a distance  $\ge 3$  splitting or for some $X$-stack $s$ and $Y$-stack $t$ and some component $R$ of $CR(s \cup t)$, $S-R$ is compressible in both $V_X$ and $V_Y$.
\endproclaim

\demo{Proof} We proceed just as with 5.1: If the splitting is not a distance $ \ge 3$ splitting we find a vertex $C \in C(S)$ disjoint from vertices $a \in K_X, b \in K_Y$ and lying in some component $R$ of $CR(s \cup t)$ for some $X$-stack $s$ and $Y$-stack $t$. If $R$ is an annulus, then $a$ and $b$ are isotopic into $S - R$; so $S-R$ is compressible in both $V_X$ and $V_Y$. \qed
\enddemo

Note : there is an algorithm for deciding whether a surface in the boundary of a handlebody is compressible in the handlebody. It seems fairly well known, but we state it for completeness. It comes from repeated applications of the following.

\proclaim{5.3 Lemma} Let $X$ be a full simplex of $C(S)$ and $F$  be a compact surface in $\partial(V_X)$ whose boundary meets $X$ efficiently. Suppose that $F$ split open along $X$ is simply connected. Then:

(1) If $F$ is compressible in $V_X$ then there is a wave of $X$ lying in $F$, and

(2) If there is a wave of $X$ in $F$, then surgery of $X$ along this wave produces a simplex $X'$ of $C(S)$ with $V_{X'} = V_X$ and $i(X', \partial(F)) \le i(X, \partial(F))-2$.
\endproclaim

By a Dehn twist along a simplex $X \subset C(S)$ we mean the product of the (commuting) Dehn twists along the vertex curves of $X$.

\proclaim{5.4 Theorem} Let $X,Y$ be simplexes of $C(S)$ with $Y$ full and having no essential vertices. Suppose that for each $X$-stack $s$, $S-s \cup X$ is simply connected. Let $h: S \to S$ be the Dehn twist along $X$. Then for $n \ge 2$
$(S;V_{h^n(Y)}, V_Y)$ is a distance $\ge 3$ splitting.
\endproclaim

\demo{Proof} We will use Theorem 5.1. The following describes $h^n(Y)$. For each component $x_i$ of $X$ let $k_i = i(x_i,Y)$ and let $A_i$ be an annular neighborhood  of $x_i$. Replace each of the $k_i$ arcs of $A_i \cap Y$ by arcs which circle $A_i$ $n$ times (and smooth to general position relative to $Y$) to get a collection of curves representing $h^n(Y)$. See Figure 5.

\bigpagebreak
\epsfysize=1.93in
\centerline{\epsffile{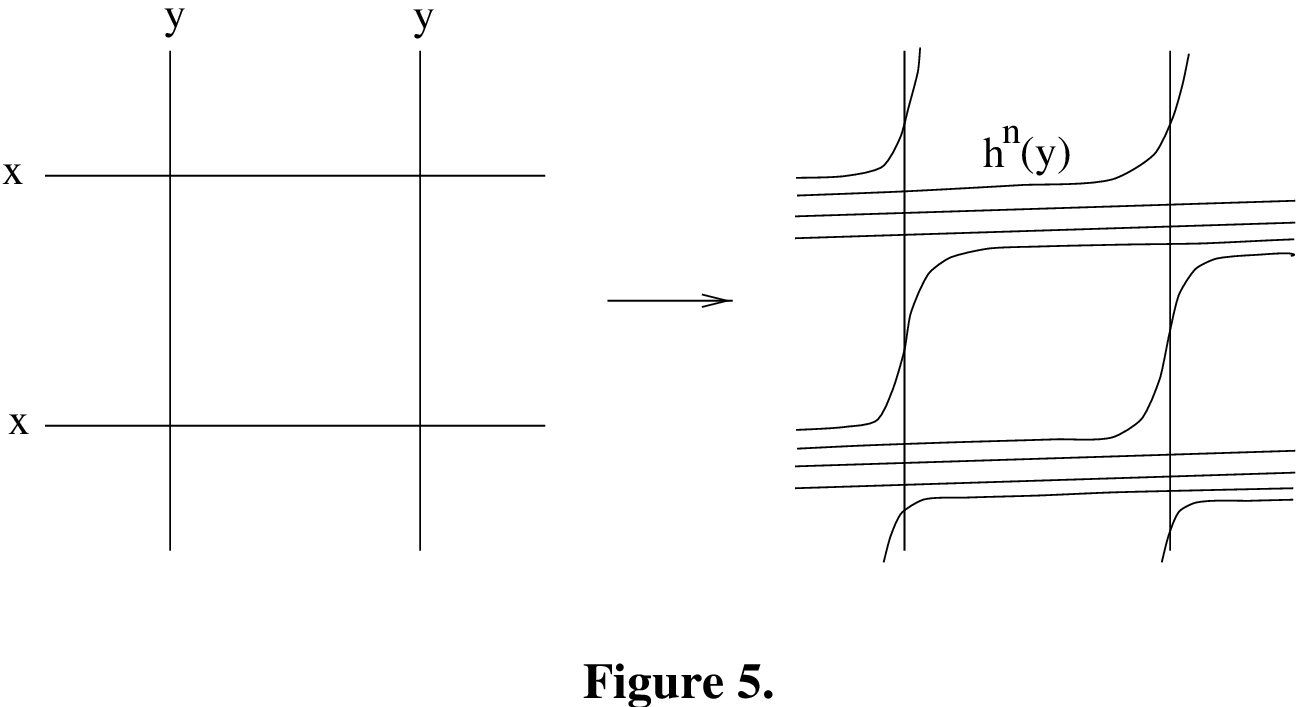}}
\bigpagebreak

The regions of $S - h^n(Y) \cup Y$ are of two types. There are the {\it old} regions which are essentially the regions of $S - X \cup Y$, but whose boundary has been twisted -- with some edges shrunk to near (old) vertices and alternate edges  expanded to near two (old) edges. Then there are the {\it new} squares. They come in {\it partial} $Y$-stacks each of which begins at an old region on one side of some $A_i$, circles $A_i$ $(nk_i -1)/k_i$. times and ends at an old region on the other side of $A_i$. There are $k_i$ of these partial $Y$-stacks in each $A_i$.

The old regions at one or both ends of a partial $Y$-stack of (new) squares may be squares; so that the (new) $Y$-stacks, relative $(h^n(Y),Y)$, will consist of these partial stacks joined together along old squares. The top and bottom will lie in old large regions. In particular, every old square in an $X$-stack, relative to $(X,Y)$ will lie in a fixed new $Y$-stack. The hypothesis necessitates that each of these $X$-stacks meets each component of $X$. So each new $Y$-stack contains partial stacks circling each $A_i$.

The $h^n(Y)$-stacks, relative to $(h^n(Y),Y)$ occupy essentially the same space as the $X$-stacks, relative to $(X,Y)$; in fact we may assume that each of the latter lies in a unique one of the former. 

So for each $h^n(Y)$-stack $s$ and each $Y$-stack $t$ , relative to $(h^n(Y),Y)$ there is an $X$-stack $s'$, relative to $(X,Y)$ and a subset of $s \cup t$ which misses any of its singular vertices and is isotopic to $s' \cup X$. Here we use the assumption $n \ge 2$ to show that each partial stack crosses $s'$ twice and so $s' \cup $ partial stack contains, up to isotopy, $s' \cup x_i$.

The conclusion now follows from Theorem 5.1 \qed
\enddemo

\proclaim{5.5 Example} The following diagram, Figure 6, satisfies the hypothesis of Theorem 5.4 and so provides  examples of  distance $\ge 3$ splittings.
\endproclaim

\bigpagebreak
\epsfysize=2.15in
\centerline{\epsffile{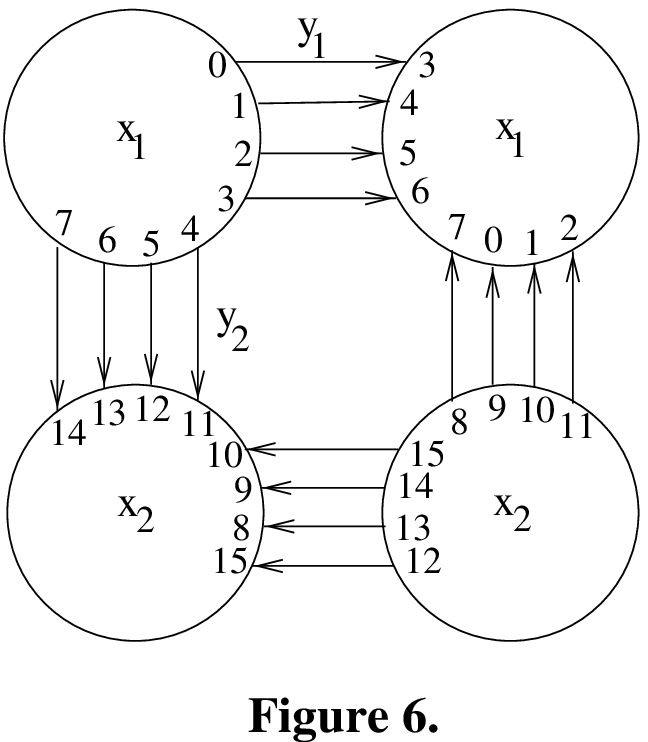}}
\bigpagebreak

{\bf 6. Distance two, genus two splittings.} In this section we describe all 3-manifolds which admit such splittings. We show that there is always a torus separating the manifold into  pieces of two specific types. Either type may reduce to a solid torus in special cases (which we describe). The details of this decomposition can be read off from a particularly nice diagram of the manifold.

One of the pieces will be Seifert fibered over $B^2$ with at most two  singular fibers. We call such a manifold a { \it generalized torus knot space} or simply a GTS and denote it by $GTS(\beta_1/\alpha_1,\beta_2/\alpha_2)$ to indicate the fiber invariants. It is the complement of an open regular neighborhood of a  ``torus knot'' (i.e. lying on a splitting torus) in a lens space. Completing $GTS(\beta_1/\alpha_1,\beta_2/\alpha_2)$ to a Seifert fibration over $S^2$ with Euler number $b$ produces the lens space
$$
L_{\alpha_1\alpha_2b-\alpha_1\beta_2 - \alpha_2\beta_1,q}.
$$
The following is immediate from the classification of Seifert manifolds.

\proclaim{6.1 Lemma} $GTS(\beta_1/\alpha_1,\beta_2/\alpha_2)$ is a solid torus if and only if $|\alpha_1|=1$ or $|\alpha_2|=1$.
\endproclaim

The other type of piece might be called a {\it one-bridge in a lens space knot complement} as it is the complement of a neighborhood of a knot which lies, except for one bridge, on a splitting torus for a lens space. By this reasoning we perhaps should be calling a GTS a {\it zero-bridge in a ...}, but we won't do either. We use the expression OBL to refer to such a manifold.

An OBL has the following structure. Let $R = T - Int(D_1 \cup D_2)$ where $T$ is a torus and $D_1$ and $D_2$ are disjoint disks in $T$. Let $a_0$ and $a_1$ be simple closed curves in $R$ which meet efficiently. Then
$$
 R \times [0,1] \cup_{a_0 \times 0} 2-handle \cup_{a_1 \times 1} 2-handle
$$
is an OBL which we denote by $OBL(a_0,a_1)$ and every OBL is so obtained. 

Note that $(T;a_0,a_1)$ is a diagram for a lens space $L$. Let $k_i$
be an arc from the center of $D_1$ to the center of $D_2$ in $T-a_i$
and let $k'_0$ be an arc in $L$ obtained by pushing $Int(k_0)$  into
$V_0$. Then $k=k'_0 \cup k_1$ is a knot in $L$ and the complement of a neighborhood of $k$ is homeomorphic to $OBL(a_0,a_1)$.

Let $<\,,\,>:H_1(S) \times H_1(S) \to {\Bbb Z}$ denote the intersection pairing on a surface $S$. So for oriented simple closed curves $x,y$ in $S$ meeting efficiently $<x,y> = i(x,y)$ means that the intersection number is $+1$ at each point of $x \cap y$.

\proclaim{6.2 Lemma} $OBL(a_0,a_1)$ is a solid torus if and only if $|<x,y>|=i(x,y)$ and either $i(a_0,a_1)=1$ or $D_1$ and $D_2$ lie in regions of $T-a_0 \cup a_1$ having an edge in common.
\endproclaim

\demo{Proof}  First we construct a genus two  diagram for $M = OBL(a_0,a_1)$. We may assume that $D_1$ and $D_2$ are (small) disks centered at points of a parallel copy $a'_0$ of $a_0$ and that the attached handle misses $a'_0 \times 0$. Then $a'_0 -Int(D_1 \cup D_2$ will be the union of two arcs $b_1$ and $b_2$ such the disks $E_i = b_i \times [0,1], \, i=1,2$ split  $V= R \times [0,1] \cup_{a_0 \times 0} 2-handle $ into a 3-cell. Thus $V$ is a handlebody and $(\partial(V);\{x_1,x_2\},y)$ is a diagram for $M$ where $x_i = \partial(E_i)$ and $y=a_1 \times 1$. See Figure 7.

\bigpagebreak
\epsfysize=1.8in
\centerline{\epsffile{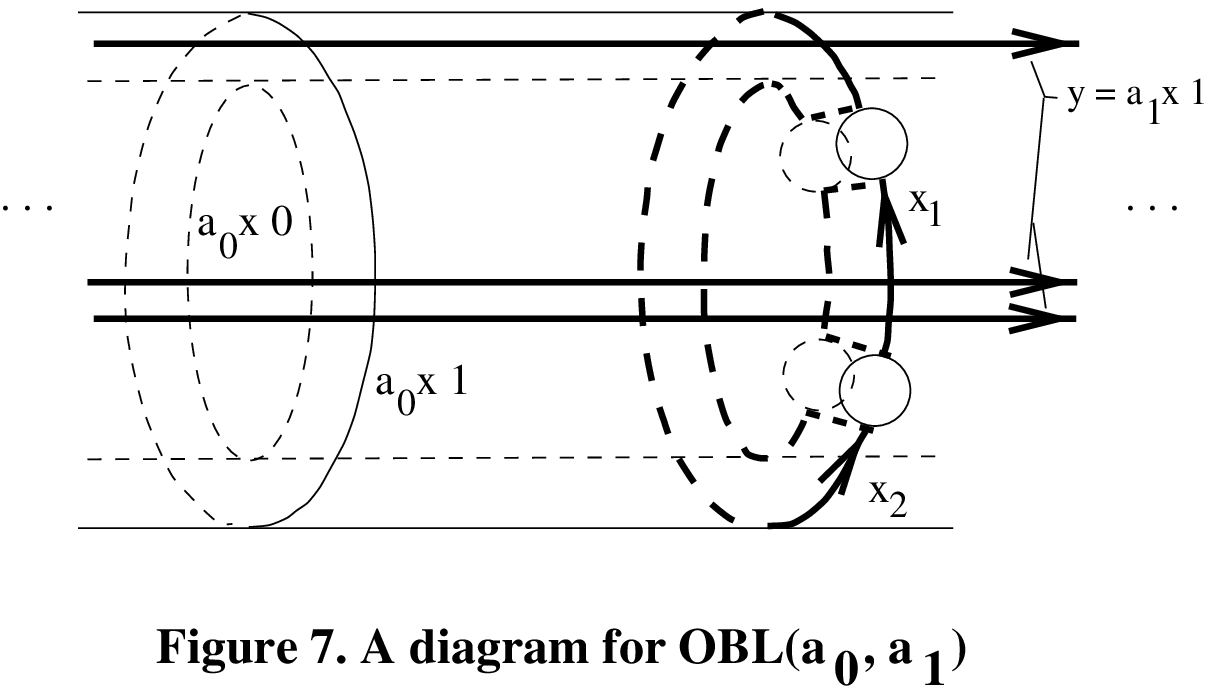}}
\bigpagebreak

Since some covering space of $M$ embeds in $S^3$, $M$ is irreducible. Thus $M$ is a solid torus if and only if $\partial(M)$ is compressible in $M$ . This in turn holds if and only if $\pi_1(M) = {\Bbb Z}$. We break the proof into two cases.

First suppose that $|<a_0,a_1>| < i(a_0,a_1)$. We will show in this case that there is no wave of $X$ in $\partial(V)-y$; so by 5.3 $\partial(V)-y$ is incompressible in $V$. It follows from a theorem of Jaco [J] (valid for a single 2-handle attachment) that $\partial(M)$ is incompressible in $M$ to complete the proof in this case.

Now geometric and algebraic intersection numbers agree for curves meeting efficiently on a torus. Thus there are exactly two bigon regions of $T-a_0 \cup a_1$ and they contain $D_1$ and $D_2$ respectively. They lie on opposite sides of both $a_0$ and  $a_1$. The situation must be as shown in Figure 8.

\bigpagebreak
\epsfysize=1.97in
\centerline{\epsffile{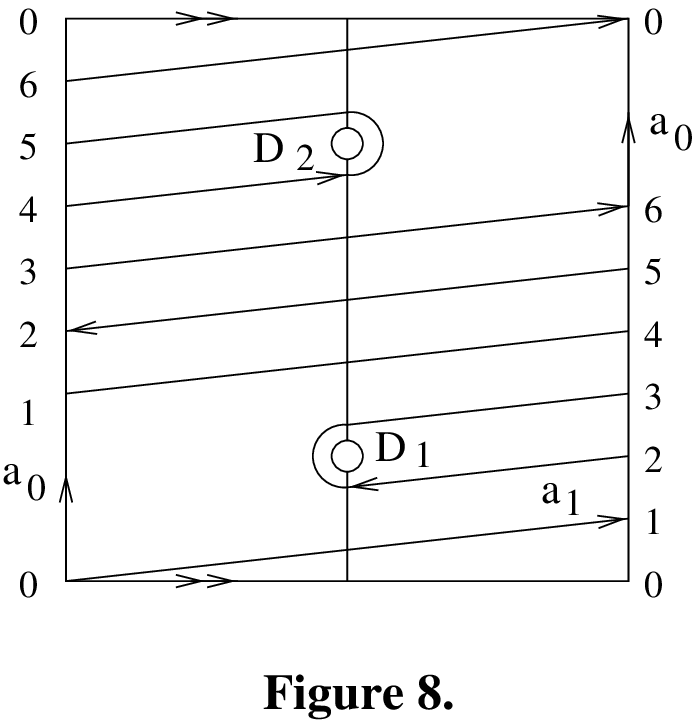}}
\bigpagebreak

Now $A= \partial(V)-R \times 1$ is an annulus. It meets each of $x_1,x_2$ in a spanning arc of $A$. The two (square) components of $A -x_1 \cup x_2$ might serve as paths for a wave to move from one part of $R \times 1$ to another. But they don't: each sees a ``dead end'' associated with one of the bigons. That is a potential wave could be isotoped to lie in $R \time 1$. But from there it can be seen not to exist at all.

Now assume $a_0$ and $a_1$ meet positively at every point. In this case there may be waves -- even when $\partial(V)-y$ is incompressible; so we need a different argument.

Let $r=i(a_0,a_1)$. If $r=0$, then $\pi_1(M) = {\Bbb Z}*{\Bbb Z}$; so $M$ is not a solid torus. If $r=1$ then $M$ is a solid torus. Thus we assume $r \ge 2$.
If $D_1$ and $D_2$ lie in the same component of $T - a_0 \cup a_1$, then $\pi_1(M)= {\Bbb Z}*{\Bbb Z}_r$; so we assume this does not happen.

Now we may assume that $T= {\Bbb R}^2/{\Bbb Z}^2$, $a_0$ is the image of the $y$-axis, $a_1$ is the image of $y=sx/r$ (where $(r,s)=1$), and $D_i$ is a small disk centered just to the right of $(0,1/2r)$ for $i=1$ and $(0,(t+1/2)r)$ for $i=2$. See Figure 9. We refer to $M$ as $OBL^+(s/r,t)$

\bigpagebreak
\epsfysize=1.7in
\centerline{\epsffile{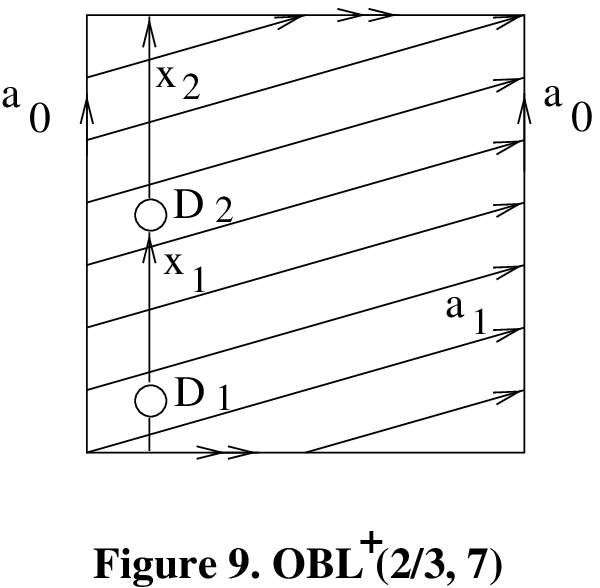}}
\bigpagebreak

Corresponding to the splitting described above we have the presentation:

$$
\pi_1(OBL^+(s/r,t))=<x_1,x_2:x_{i_1}x_{i_2} \dots x_{i_r} = 1>
$$
where $i_j = 1$ if $[1+(j-1)s] \in \{1,2, \dots ,t\}$ and $i_j = 2$ otherwise.
Here $[n]$ denotes the least non-negative residue of $n$ modulo $r$. We will complete the proof in this case by showing that $\pi_1(M) \cong {\Bbb Z}$ if and only if $t \equiv \pm 1,\pm s \text{ mod } r$.

If $t= \pm 1$ then the relation involves just one occurrence of $x_1$ or $x_2$, and clearly $\pi_1(M)={\Bbb Z}$. The case $t = \pm s$ reduces to the case $t = \pm 1$ on interchanging the roles of $a_0$ and $a_1$.

By obvious symmetries (which may reverse orientation) we may assume that
$1 < t \le r/2$ and $1\le s < r/2$. By combining terms, and conjugating, if necessary, we can express the relation in the form
$$
x_1^{n_1}x_2^{m_1}x_1^{n_2}x_2^{m_2} \dots x_1^{n_k}x_2^{m_k}
$$
where $n_i, m_i \ge 1$. The $n_i$ (respectively $m_j$) are the lengths of maximal sequences $[q], [q+s], [q+2s], \dots$ lying in $\{1,2,\dots,t\}$ (respectively $\{t+1,t+2,\dots, r\}$).

If $s < t$ (and so $s<r-t$) then some $n_i \ge 2$ and some $m_j \ge 2$. If $s > t$, then, say, every $n_i =1$, but using $s<r/2$ we can show that some two $m_j$'s differ by at least two. In either case Lemma 6.3 below applies to show
$\pi_1(M) \not\cong {\Bbb Z}$. \qed
\enddemo

\proclaim{6.3 Lemma} Let $G=<x_1,x_2: x_1^{n_1}x_2^{m_1}x_1^{n_2}x_2^{m_2} \dots x_1^{n_k}x_2^{m_k}=1>$ be a group presentation with $n_i >0, m_j >0$ for all $i,j$.
If either some $n_i \ge 2$ and some $m_j \ge 2$ or, say, some two $m_j$'s differ by at least $2$ , then the Alexander polynomial $\Delta_G(z)$ is not constant. So $G \not\cong {\Bbb Z}*H$ for any group $H$.
\endproclaim

\demo{Proof} Let $n= \sum n_i,\, m=\sum m_i$, $d= (n,m)$, $n=d\bar n$, and $m = d\bar m$. Then $G/G' = {\Bbb Z} \times {\Bbb Z}_d$ and we have an epimorphism
$\phi:G \to <z:\,\,>$ of $G$ to the (multiplicative) infinite cyclic group given by $\phi(x_1)=z^{\bar m}, \phi(x_2)=z^{-\bar n}$ which is unique up to the automorphism $ z \to z^{-1}$.
We use the same notation for the extension ${\Bbb Z}G \to {\Bbb Z}[z,z^{-1}]$.

The Alexander polynomial, $\Delta_G(z)$, of $G$ is a generator for the ideal
$$
(\phi(\partial r/\partial x_1),\phi(\partial r/\partial x_2).
$$

Now $ (\partial r/\partial x_1)(x_1-1)= (\partial r/\partial x_2)(1-x_2)$.
Since  $(\bar n, \bar m) =1$, the greatest common divisor of
$$
\phi(x_1 -1) = z^{\bar m}-1 \text{ and } \phi(1-x_2) =1- z^{-\bar n}  = z^{-\bar n}(z^{\bar n} -1)
$$
is $z-1$. It follows that
$$
\Delta_G(z)=\phi(\partial r/\partial x_1)/(1+z+ \dots +z^{\bar n -1}) =
\phi(\partial r/\partial x_2)/(1 + z+ \dots +z^{\bar m -1}).
$$

Now for the first part of the theorem we may assume that $\bar n \ge \bar m$ and that some $m_j \ge 2$. After cyclically permuting, if necessary, we may further assume $r = wx_1x_2^{m_j}x_1u$ for some positive words $w,u$ in the generators.

Then $\partial r/\partial x_1$ contains the terms $w$ and $wx_1x_2^{m_j}$ which, since all exponents are positive, are not cancelled by any other terms.
The difference of degrees of the images of these terms under $\phi$ is
$|\bar m - \bar n m_j| \ge \bar n$. Thus we are left with a non-constant (Laurant) polynomial after dividing by $1+z+ \dots +z^{\bar n-1}$.

For the second part suppose $n_1=n_2= \dots =1$,  that $m_1 $ is the smallest
$m_j$, and that $m_j-m_1 \ge 2$ for some $j$. We replace $x_1$ by $x = x_1x_2^{m_1}$ to get a new presentation  $G = <x,x_2:x^2x_2^{m_2-m_1} \dots xx_2^{m_j-m_1} \dots>$ to which the first part applies.

For the final part suppose $G \cong {\Bbb Z}*H$ for some group $H$. Then $H/H' \cong {\Bbb Z}_d$. So elements of $H$ map to $1$ under $\phi:G \to <z:\,\,>$.
We get another presentation for $G$ by adding a generator, but no relations, to a presentation for $H$. The entries of the Jacobian matrix of this presentation map to integers under $\phi$. Thus $\Delta_G(z)$ is constant. \qed

\enddemo

\proclaim{6.4 Theorem} Let $(S;V_1,V_2)$ be a genus two, distance two splitting of a 3-manifold $M$. Then there is a torus $T \subset M$ splitting $M$ into a GTS, $M_1$, and an OBL, $M_2$. There is a diagram $(S;\{x_1,x_2\},\{y_1,y_2\})$ for the splitting so that $S-x_1 \cup y_1$ contains an essential, nonseparating simple closed curve $z$, and for any such diagram
$$
M_1 = GTS(\beta_1/<z,x_2>,\beta_2/<z,y_2>)
$$
and
$$
M_2 = OBL(S-z;x_1,y_1)
$$
\endproclaim

\demo{Proof} Let $(S;X,Y)$ be any diagram for the splitting. By assumption there are vertices $x \in K_X, y\in K_Y$ with $d(x,y) = 1$. So there is an essential simple closed curve $z \subset S-x\cup y$. We may assume $z$ does not separate $S$; for otherwise $x \cup y$, being connected, would lie in one component of $S-z$ and we could replace $z$ by an essential, nonseparating simple closed curve in the other component. Similarly we may assume $x$ and $y$ are nonseparating. For if, say,  $x$ separates $S$ then each component of $S-x$ contains a nonseparating vertex of $K_X$ and we could replace $x$ by one missing $z$.
Thus we can extend $x,y$ to a diagram as desired.

Now given any such diagram take disks $D_i \subset V_1, E_j \subset V_2$ bounded respectively by $x_i, y_j$. Let $A$ be a regular neighborhood of $z$ in $S$, and let $R=Cl(S-A)$ -- a twice punctured torus. A regular neighborhood $M_2$ of $R \cup D_1 \cup E_1$ in $M$ is an $OBL(x_1,y_1)$.

$M_1 = Cl(M-M_2)$ is the union of two solid tori: $V_1$ split along $D_1$ and $V_2$ split along $E_1$ pushed slightly away from $R$. These solid tori meet along $A$ which circles them $<z,x_2>$ and $<z,y_2>$ times respectively. 

Of course, the $\beta_i$'s can also be read off from the diagram as can the gluing map $\partial(M_1) \to \partial(M_2)$. Perhaps it is better not to try to squeeze too much in the statement of the theorem, but to leave it to calculation as illustrated in the examples below. \qed
\enddemo

\proclaim{6.5 Comments} 

1. An OBL may be Seifert fibered. The ``one bridge'' might be isotoped back onto the torus $T$. One condition for being able to do this for $OBL^+(s/r,t)$ is:
one of $A_i \cap B_j = \emptyset$ where
$$
A_1=\{1,2,\dots,t-1\}, A_2=\{t+1,t+2,\dots,r-1\}
$$
$$
B_1=\{s,2s,\dots,(\bar s-1)s\}, B_2=\{(\bar s +1)s, (\bar s +2)s, \dots,(r-1)s\};\,\bar s s \equiv t \text { mod } r.
$$
In particular every $OBL^+(1/r,t) $ is Seifert fibered, as is, for example $OBL^+(3/10,6)$.

2. The complement (in $S^3$) of a 2-bridge knot is $OBL(a_0,a_1)$ where there is a simple closed curve $a'_1$ isotopic to $a_1$ in $T$ such that $i(a_0,a'_1) =1$ and $a'_1$ meets the projection of the knot in $T$ in a single point.

3. There are OBL's with $\Delta(0) \ne \pm 1$ and which therefore cannot be surface bundles over $S^1$; for example $OBL^+(3/14,7)$.
\endproclaim

\proclaim{6.6 Examples} Each of the diagrams of Figure 10 satisfies the conditions of Theorem 6.4.

In (a) there is a torus splitting the manifold into a $GTS(1/2,1/3)$ and an $OBL^+(1/5,3)$ which by 6.5(1) is also Seifert fibered. This is a non-trivial graph manifold.

In (b) we have a splitting into a $GTS(1/3,1/4)$ and an $OBL^+(1/3,1)$ which by 6.2 is a solid torus. Thus the manifold is Seifert fibered. This could also be seen from Theorem 3.2 as the natural way of completing the diagram to pants decompositions will give $i(x_3,y_3)=2$.

In (c) we have a splitting into a $GTS(1/1,1/3)$, which is a solid torus, and an $OBL(x_1,y_1)$, which by 6.5(2) is a 2-bridge knot complement. With a little more effort this can be seen to be $2/3$ surgery on the figure eight knot.
\endproclaim

\epsfysize=5.03in
\centerline{\epsffile{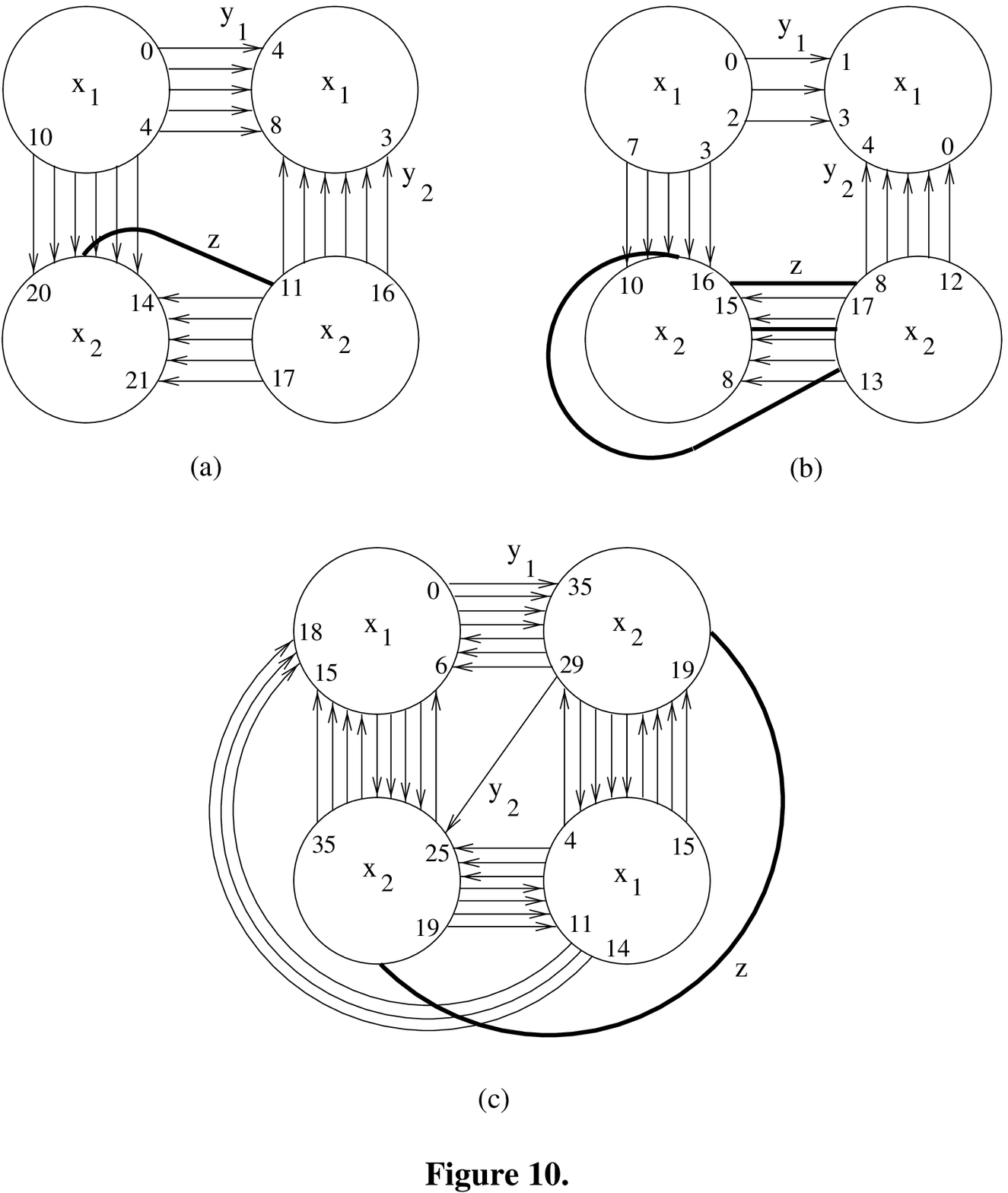}}
\bigpagebreak

\enddocument